\documentstyle[12pt]{amsart}
\newcommand{\der}{\partial}

\newtheorem{thm}{Theorem}[section]
\newtheorem{crl}[thm]{Corollary}
\newtheorem{prp}[thm]{Proposition}
\newtheorem{lm}[thm]{Lemma}

\newcommand{\D}{\mathcal D}

\begin{document}
\title[Identities and derivations for jacobians] 
{Identities and derivations for Jacobian algebras}
\author{A.S. Dzhumadil'daev}
\address{International Centre for Theoretical Physics, Trieste, Italy, \newline
Institute of Mathematics, Academy of Sciences of Kazakhstan, {\hbox{Almaty}}}
\email{askar@@math.kz}
\maketitle

\begin{abstract} Constructions of $n$-Lie algebras by strong $n$-Lie-Poisson algebras are given. First cohomology groups of adjoint module of  Jacobian 
algebras are calculated. Minimal identities of $3$-Jacobian algebra are found. 
\end{abstract}
 
\section{Introduction} 
Let $U$ be an associative commutative algebra over a field $K$ 
with commuting derivations 
$\der_1,\ldots,\der_n.$ Say, $U$ is an algebra of differentiable functions 
on $n$-dimensional manifold or polynomial algebra $K_n^+=K[x_1,\ldots,x_n]$
or Laurent polynomial algebra $K_n=K[x_1^{\pm 1},\ldots,x_n^{\pm 1}].$ 
In these examples, $\der_i=\der/\der\,x_i$ are partial derivations. 
If not stated otherwise, the characteristic $p$ of the field $K$ 
is supposed to be $0.$ 

Let $Jac_n^S: \wedge^n U\rightarrow U$ be Jacobian map: 
$$Jac_n^S(u_1,\ldots,u_n)= det\,(\der_iu_j)=\left|\begin{array}{ccc}
\der_1u_1&\cdots&\der_1u_n\\ \hdots&\hdots&\hdots\\ \der_nu_1&\cdots&\der_nu_n\\
\end{array}\right|$$
Define $(n+1)$-linear map $Jac_{n+1}^W:\wedge^{n+1}U\rightarrow U$ by 
$$Jac_{n+1}^W(u_0,u_1,\ldots,u_n)=\left|\begin{array}{cccc}
u_0&u_1&\cdots&u_n\\
\der_1u_0&\der_1u_1&\cdots&\der_1u_n\\ \hdots&\hdots&\hdots&\hdots\\ 
\der_nu_0&\der_nu_1&\cdots&\der_nu_n\\
\end{array}\right|$$
In terms of wedge products we see that 
$$Jac_n^S=\der_1\wedge\cdots\wedge \der_n,$$
$$Jac_{n+1}^W=id\wedge\der_1\cdots\wedge \der_n,$$
where $id:U\rightarrow U, u\mapsto u,$ is the identity map. 

In \cite{Filipov1}, \cite{Filipov2}, \cite{DalTakh} are proved that the $n$-ary multiplication 
$\omega=Jac_n^S$ satisfies the  identity 
\begin{equation}\label{FI1}
g^{\omega}=0,
\end{equation}
where
$$g^{\omega}(u_1,\ldots,u_{2n-1})=
\omega(u_1,\ldots,u_{n-1},\omega(u_n,\ldots,u_{2n-1}))-$$
$$
\sum_{i=1}^n\omega(u_1,\ldots,u_{i-1},\omega(u_1,\ldots,u_{n-1},u_i),
u_{i+1},\ldots,u_{2n-1}).$$
In \cite{Takh} this identity is called {\it a Fundamental Identity}. 
We call it {\it Fundamental Identity of type {\rm I}.}   

We notice that the Jacobian $\omega=Jac_n^S$ 
satisfies one more identity, that we call 
{\it a Fundamental identity of type {\rm II:}}
\begin{equation}\label{FI2}
f^{\omega}=0,
\end{equation}
where
$$f^{\omega}(u_1,\ldots,u_{n-1},v_1,\ldots,v_{n+1})=$$
$$
\sum_{i=1}^{n+1}(-1)^i\omega(u_1,\ldots, u_{n-1},v_i)\cdot\omega(v_1,\ldots,\hat{v_i},\ldots,v_{n+1}).
$$

It is known that, the Jacobian $\omega=Jac_n^S$  satisfies the Leibniz 
identity for the multiplication $\cdot:$
\begin{equation}
\label{Leibniz}
\omega(u\cdot u',u_2,\ldots,u_n)=u\cdot \omega(u',u_2,\ldots,u_n)+
u'\cdot \omega(u,u_2,\ldots,u_n).
\end{equation}
Here $u,u',u_1,\ldots,u_{2n+1}, v_1,\ldots,v_{n+1}$ are any elements of $U.$ 

In our paper we consider algebras with many operations. 
An operation or multiplication on algebra is a polylinear map. 
If $\omega: U\times\cdots\times U\rightarrow U$ is a polylinear map with 
$n$ arguments, then $\omega$ is called a {\it $n$-ary multiplication} 
on $U.$ The space of $n$-ary polylinear maps on $U$ is denoted by $T^n(U,U).$ 
If $n=0,$ then we set $T^0(U,U)=U.$ The set of operations on $U$ is called  
{\it a signature} of algebra \cite{Kurosh}.
The algebra with vector space $U$ and signature 
$\Omega=\{\omega,\eta,\ldots\}$ is denoted as 
$(U,\Omega)$ or $(U,\omega,\eta,\ldots)$ or just $U,$ 
when it is clear which multiplications are considered. 

An algebra $(U,\omega)$ with skew-symmetric $n$-ary multiplication 
$\omega$ that satisfies (\ref{FI1}) is called {\it $n$-Lie} \cite{Filipov1}. 
An algebra $(U,\cdot,\omega)$ is called {\it $n$-Lie-Poisson,} if 
$(U,\cdot)$ is an associative commutative algebra and it satisfies 
the identities (\ref{FI1}) and (\ref{Leibniz}). If it satisfies one more identity, namely the  identity (\ref{FI2}), then this algebra is called {\it strong $n$-Lie-Poisson.} 

Sometimes $n$-Lie algebras are called  Nambu \cite{Nambu}, Filipov or 
Takhtajan algebras and $n$-Lie-Poisson algebras are called  
Nambu-Poisson algebras. 

{\bf Question.} Does the fundamental identity of type {\rm II} 
follow  from the fundamental identity of type {\rm I}
and from the Leibniz identity, if $n>2$ and $p=0$ ? 

Here, we suppose that commutativity and associativity identities for the 
binary multiplication and skew-symmetric identity for $n$-multiplication 
are given. Otherwise, more likely,  the answer 
to this question will be negative. 

For some statement $\mathcal X$ denote by $\delta(\mathcal X)$ its 
Kroneker symbol: 
$\delta({\mathcal X})=1,$ if $\mathcal X$ is true and $\delta({\mathcal X})
=0,$ if $\mathcal X$ is false. 
For some set of vectors $Y$ denote by $<Y>$ its  linear span. 

Let  
$$r^{\omega}(u_1,\ldots,u_{n-2},u_{n-1},u_n,\ldots,u_{2n})=$$
$$\sum_{i=n-1,n+1,\ldots,2n}(-1)^{i+n+\delta(i\le n-1)}\omega(u_1,u_2,\ldots,u_{n-2},u_i)\omega(u_n,u_{n+1},\ldots,\hat{u_i},\ldots,u_{2n})+$$
$$\sum_{i=n-1,n+1,\ldots,2n}(-1)^{i+n+\delta(i\le n-1)}\omega(u_n,u_2,\ldots,u_{n-2},u_i)\omega(u_1,u_{n+1},\ldots,\hat{u_i},\ldots,u_{2n}).$$
So, $r^{\omega}$ is the polynomial with $n$ skew-symmetric arguments $u_{n-1},u_{n+1},\ldots,u_{2n}$ and two symmetric arguments $u_1$ and $u_n.$ 

By theorem 6 \cite{DalTakh} the answer to our question will be positive, 
if the inverse of this statement is true: 
the identity $f^{\omega}=0$ follows from the identity $r^{\omega}=0$ and 
Leibniz rule (\ref{Leibniz}).    

We give a positive answer to the last statement for $n=3$ 
(theorem~\ref{minimal}). 
If $p>0,$ it is more likely that the answer would be negative. 
For $n=2$ the answer is also negative. There exist examples of Lie-Poisson algebras, that are not strong. For example, $(K[x_1,x_2,x_3,x_4],\partial_1\wedge \partial_2+\partial_3\wedge\partial_4)$ is  such an algebra. 
 
Let  $(A,\Omega)$ be an algebra with some vector space 
$A$ and signature $\Omega.$ For $\Omega'\subseteq \Omega$ 
a linear map $D: A\rightarrow A$ is called $\Omega'$-derivation, 
if 
$$D(\omega(a_1,\ldots,a_n))=\sum_{i=1}^n\omega(a_1,\ldots,a_{i-1},D(a_i),
\ldots,a_n),$$ 
for any $a_1,\ldots,a_n\in A$ and for all $\omega\in \Omega'.$ 
Call $\Omega$-derivation a derivation.  
Let $Der\,(A,\Omega')$ be a space of all derivations of $(A,\Omega').$ 
Set $Der\,A= Der\,(A,\Omega).$  

In terms of operators $L_{a_1,\ldots,a_{n-1}}: A\rightarrow A, 
a\mapsto \omega(a_1,\ldots,a_{n-1},a)$ we see that $\omega$ satisfies (\ref{FI1}), if 
and only if 
$$L_{a_1,\ldots,a_{n-1}}\in Der A,$$
for any $a_1,\ldots,a_{n-1}\in A.$ Derivations of the form $L_{a_1,\ldots,a_{n-1}}$ 
are called {\it interior} derivations. Let $Int\,A$ be a space of interior derivations. If $A$ is $n$-Lie, then $Int\,A$ is a Lie algebra under commutator 
of operators. 
Moreover, $Int\,A$ is an ideal of $Der\,A:$ 
$$[D,L_{a_1,\ldots,a_{n-1}}]=\sum_{i=1}^nL_{a_1,\ldots,a_{i-1},D(a_i),a_{i+1},\ldots,a_n},$$
for any $D\in Der\,A, a_1,\ldots,a_n\in A.$ In particular one can consider 
a factor-algebra, an algebra of outer derivations, $Out\,A=Der\,A/Int\,A.$  

In \cite{Filipov1} two examples of $n$-Lie algebras were constructed. 
The first one is vector products algebra and the second one is 
Jacobian algebra. In this paper it was established that any 
derivation of vector products algebra  is interior.
In \cite{Filipov1} it was also noticed that, if $(A,\omega)$ is $(n+1)$-Lie, 
then $(A,i(a)\omega)$ is $n$-Lie, where $n$-ary map $i(a)\omega$ is defined by 
$$i(a)\omega(a_1,\ldots,a_n)=\omega(a,a_1,\ldots,a_n).$$

In our paper we give a generalisation of Jacobian algebras. Namely, 
we establish that the algebra $(U,Jac_{n+1}^W)$ becomes a 
$(n+~1)$-Lie algebra. If $U$ has unit element $1,$ then $\der_i(1)=1,$ for any 
$i=1,\ldots,n.$ If $U$ has unit, then 
the $(n+1)$-ary algebra $(U,Jac_{n+1}^W)$ allows us to obtain  
the $n$-ary Jacobian algebra 
$(U,Jac_n^S)$ by the restriction operation: $Jac_n^S=i(1)Jac_{n+1}^W.$

The theory of polynomial identities is well developed for ordinary algebras, 
i.e., for algebras with binary operations. The case of multi operation  
algebras needs some detailed information about $n$- or  $\Omega$-words and 
$\Omega$-polynomials. Necessary 
definitions and descriptions 
of $\Omega$-words are given in section \ref{oneone}
(theorem \ref{22may}). We introduce a notion of $\Omega$-degree for  
$\Omega$-polynomial, that is, the number of operations. For example, Jacobian 
algebra considered as an $n$-Lie algebra, has only one $n$-ary operation, 
denoted by $\mu_n=Jac_n^S,$ and it has only one identity of $\Omega
\mbox{-degree}$ 
1 (skew-symmetric identity for Jacobian) and only one identity of 
$\Omega$-degree 2 ($n$-Lie identity). The identity $f^{\mu_n}=0$ is not 
an $n$-Lie identity, since the construction of $f^{\omega}$ needs  
binary  operation. 
Probably Jacobian algebra as $n$-Lie algebra has no any identity of 
$\Omega$-degree 3. 
If we consider Jacobian algebras as $n$-Lie-Poisson algebras, i.e., as 
algebras with one binary operation $\mu_2: (u,v)\mapsto u\cdot v$ and one 
$n$-ary operation $\mu_n=Jac_n^S,$ then it has two identites of 
$\Omega$-degree 1 
(commutativity for $\mu_2$ and skew-symmetry for $\mu_n$) 
and three identities of $\Omega$-degree 2
($n$-Lie for $\mu_n,$ Leibniz rule between $\mu_2$ and $\mu_n,$ and 
associativity identity for $\mu_2$). As $n$-Lie-Poisson algebra,  
Jacobian algebra has one more identity of $\Omega$-degree 3
(identity between two $\mu_n$-s and one $\mu_2$). It seems that both 
directions of studying $\Omega$-identities of Jacobian algebras will be 
very interesting. In our paper we describe $\Omega$-degree 2 identities for 
the Jacobian algebra $(K_3,Jac_3^S)$ as $3$-Lie algebra.

In cocycle and identity constructions we use two methods: 
a polynomial principle and $\D$-invariants method (section \ref{poly}).  
The polynomial principle allows us to restrict our considerations by 
the case of polynomial algebras. As it turns out, almost all of our 
cocycles and polynomials are $\D$-invariant. To prove that  
$\D$-invariant polynomial is an identity it is sufficient to calculate 
these polynomials on their supports. In other words, the 
identity (cocyclicity) checking problem we reduce to the calculation 
problem whether a polynomial is equal to 0 on some concrete arguments. 
Here the use of computer calculations are very helpfull.

Let $U=K_n.$ In section~\ref{SS} 
we prove that the classes of the  following linear maps 
consist of the basis of $Out\,(U,Jac_n^S):$
$$\Delta:=\sum_{i=1}^nx_i\der_i+n(1-n)^{-1},$$
$$D_{-\theta}: x^{\alpha}\mapsto \delta_{\alpha,-\theta},$$ 
$$D_i:=x^{-\theta+\epsilon_i}\der_i,\quad i=1,\ldots,n,$$
where $\theta=(1,\ldots,1)$ and $\epsilon_i=(0,\ldots,\mathop{1}\limits_i,
\ldots,0).$  
In section \ref{WW} we establish that all derivations of $(U,Jac_{n+1}^W)$ 
are interior. In these sections we also prove that the algebras of interior 
derivations are isomorphic to Cartan Lie algebras of types $W$ and $S,$ namely,
  $Int\,(U,Jac_{n+1}^W)\cong W_n(U)$ and  
$Int\,(U,Jac_{n}^S)\cong S_n(U).$ These isomorphims explain index 
notations on $Jac_{n+1}^W$ and $Jac_n^S.$
 
In terms of $n$-Lie cohomology \cite{Cautheron}, \cite{Takh1}, in sections 
\ref{SS} and \ref{WW} we calculate first cohomology groups of Jacobian 
algebras with coefficients in adjoint module. Derivations of Lie algebras 
appear in a natural way in considering central extensions of Lie algebras. 
Derivations  of Lie algebras $H_1$ and $W_n$ in this sense were described 
in \cite{DzhCentral}.  Our results in the case $n=2$ are compatible with 
the results of this paper.

\section{$\Omega$-words \label{oneone} }
Let $\bf Z$ be the ring of integers, 
${\bf Z}_+=\{i\in {\bf Z}: i\le 0\}$ and 
${\bf Z}^{+}=\{i\in {\bf Z}: i>0\}.$ 
Let us given some alphabet $\aleph$ with a map 
$\aleph\rightarrow {\bf Z}_+, \alpha\mapsto |\alpha|,$ 
called  {\it arity} map. Let 
$${\frak X}=\{\alpha\in \aleph: |\alpha|=0\},$$
$$\Omega=\{\omega\in \aleph: |\alpha|>0\}.$$
Thus, $\aleph=\Omega\cup {\frak X}.$ Denote elements of $\Omega$ by  
$\omega_1,\omega_2,\ldots$ and elements of ${\frak 
X}$ by $x_1,x_2,\ldots.$

Define a {\it weight map}  
$$\aleph\rightarrow {\bf Z}, \quad \alpha\mapsto ||\alpha||,$$
by 
$$||\alpha||=1-|\alpha|.$$

\begin{lm} \label{1sept} The following conditions are equivalent
\begin{itemize}
\item $||\alpha||\ge 1$
\item $||\alpha||=1$
\item $\alpha\in {\frak X}$
\item $|\alpha|=0$
\end{itemize} 
\end{lm}

{\bf Proof.} Evident. 

Let $\Gamma(\aleph)=\{a=\alpha_1\alpha_2\cdots\alpha_k \}$ be the
set of sequences of elements of $\aleph.$ 
If $a=\alpha_1\ldots\alpha_k, b=\beta_1\ldots\beta_s\in\Gamma(\aleph),$ 
then, by definition, $a=b$ if and only if $k=s$ and $\alpha_1=\beta_1,\ldots,\alpha_k=\beta_k.$  Call the number of
elements $k$ of the sequence $\alpha=\alpha_1\cdots\alpha_k\in
\Gamma(\aleph)$ a {\it length} of $\alpha$ and denote it by
$l(\alpha).$ Prolong the weight map $||\;||$ to $$||\;||:
\Gamma(\aleph)\rightarrow {\bf Z},$$ by $$||\alpha_1\cdots
\alpha_k||=||\alpha_1||+\cdots+||\alpha_k||.$$

Let ${\bf Z}^{\infty}=\{(i_1,i_2,\ldots ): i_1,i_2,\ldots \in {\bf Z}\}.$ 
Define a map $\mu: \Gamma(\aleph)\rightarrow {\bf Z}^{\infty}$ 
by 
$$\mu(a)=
(\mu_k(a),\mu_{k-1}(a),\ldots,\mu_1(a)),$$
$$\mu_i(a)=||\alpha_i||+\cdots+||\alpha_k||, \quad i=1,2,\ldots,k,$$
if $a=\alpha_1\ldots \alpha_k, \alpha_i\in \aleph, i=1,2,\ldots,k.$

{\bf Definition.} Let $\Gamma_1(\aleph)$ be the subset of $\Gamma(\aleph),$ 
that consists of elements $a\in \Gamma(\aleph),$ such that
\begin{itemize}
\item $||a||=1$
\item if $a=\alpha_1\cdots \alpha_k,$ $\alpha_i\in\aleph, i=1,2,\ldots,k,$
then $\mu_i(a)\ge 1,$ for any $i=1,2,\ldots,k.$
\end{itemize}

{\bf Example.} Let $\Omega=\{\omega_3,\omega_2,\omega_2':
|\omega_3|=3, |\omega_2|=|\omega_2'|=2.$ 
Let 
$a=\omega_3\omega_2x_1x_2x_3\omega_2'x_4x_5\in \Gamma(\aleph),$
$b=\omega_3\omega_2x_1x_2\omega_3x_3x_4x_5\in \Gamma(\aleph).$
Then $\mu(a)=(1,2,1,2,3,4,3,1)$ and $\mu(b)=(1,2,3,1,2,3,2,0).$
Therefore, $a\in\Gamma_1(\aleph)$ and $b\not\in\Gamma_1(\aleph).$ 

\begin{lm}\label{1septe}
Elements of $\Gamma_1(\aleph)$ have the
following properties
\begin{itemize}
\item $||\alpha||=1, \alpha\in \Gamma_1(\aleph)
\Rightarrow \alpha\in {\frak X},$
\item $\Omega\not\subset \Gamma_1(\aleph),$
\item any element of
$\Gamma_1(\aleph)$ with length more than 1 begins with some element of
$\Omega.$
\item any element of $\Gamma_1(\aleph)$ ends by some element of $\frak
X.$
\end{itemize}
\end{lm}

{\bf Proof.} 
Let $a=\alpha_1\ldots\alpha_k$ be the element of $\Gamma_1(\aleph)$ with 
length $k=l(a).$ 

By definition,  
$$\mu_k(a)=||\alpha_k||\ge 1$$
Thus, by lemma \ref{1sept}, $\alpha_k\in {\frak X}.$ So, we 
have proved that any element of $\Gamma_1(\aleph)$ ends by element 
of $\frak X.$ In particular, any element of $\Gamma_1(\aleph)$ with length 1 
is an element of $\frak X.$  

Suppose that $l(a)=k>1.$ By definition,
$$\mu_1(a)=||\alpha_1\ldots\alpha_k||=1,\quad ||\alpha_2\ldots\alpha_k||\ge 1.$$
Therefore,
$||\alpha_1||\le 0.$
In other words, $\alpha_1\in \Omega.$ So, we established that any element of 
$\Gamma_1(\aleph)$ with length $>1$ begins with element of $\Omega.$ 

Define a set of $\Omega$-words \cite{Kurosh}.

{\bf Definition.} \begin{enumerate}
\renewcommand{\theenumi}{\roman{enumi}}
\item Any element of $\frak X$ is an $\Omega$-word.
\item If $a_1,\ldots,a_{k}$ are $\Omega$-words,
then $\omega a_1\ldots a_{k},$ where $|\omega|=k,$ is also a
$\Omega$-word.
\item Any $\Omega$-word is obtained by these two rules.
\end{enumerate}

Let $a=\omega a_1\ldots a_k$ be some word and $\alpha$ is a word
or element of $\Omega.$ We say that $\alpha$ enter to $a$ or that
$\alpha$ is a part of the word $a$ and write $\alpha\in a,$  if
one of the following cases are fulfilled,
\begin{itemize}
\item $\alpha$ is a word and $\alpha=a,$
\item $\alpha$ is a word and $\alpha$ is a part of $a_s$ for some
$s=1,\ldots,k,$
\item $\alpha\in \Omega$ and $\alpha=\omega,$
\item $\alpha\in \Omega$ and $\alpha$ is a part of $a_s$ for some
$s=1,\ldots,k.$
\end{itemize}

For the word $a=\omega\,a_1,\ldots,a_{k}$ define $\omega deg\,a$
or {\it $\Omega$-degree} of $a,$ by
\begin{itemize}
\item $\omega deg\,a=\sum_{j=1}^{k}\omega deg\,a_j+1,$
\item $\omega deg\,x=0, \quad x\in {\frak X},$
\item $\omega deg\,\omega=1, \quad \omega\in \Omega.$
\end{itemize}
So, $\Omega$-degree of $a$ is the number of elements of $\Omega$
that enter to $a:$ $$\omega deg\,a=|\{\omega\in \Omega\cap a\}|.$$

Let $a=\omega a_1\ldots a_k$ be some word. Call $xdeg\,a$ or {\it
$\frak X$-degree} of $a$ the number of elements of $\frak X$ that
enter to $a:$
\begin{itemize}
\item $xdeg\,a=\sum_{j=1}^{k}xdeg\,a_j,$
\item $xdeg\,y=1, \quad y\in {\frak X},$
\item $xdeg\,\omega=0, \quad \omega\in \Omega.$
\end{itemize}

A {\it degree} of $a$ is defined by $deg\,a=xdeg\,a+\omega deg\,a.$

\begin{thm} \label{22may}
The set of $\Omega$-words coincides with $\Gamma_1(\aleph).$
\end{thm}

{\bf Proof.} Denote by $\bar\Gamma_1(\aleph)$ the set of $\Omega$-words. 

Prove that $\bar\Gamma_1)(\aleph)\subseteq \Gamma_1(\aleph).$ 
Let $a\in\bar\Gamma_1(\aleph).$ We use induction on $l(a).$ If $l(a)=1,$ then 
$a=x\in\frak X.$ Thus, $||x||=1,$ and $x\in \Gamma_1(\aleph).$ Suppose that 
any element of $\bar\Gamma_1(\aleph)$ with length $<l(a)$ belongs  to $\Gamma_1(\aleph).$ 
If $a=\omega a_1\ldots a_k,$ then $l(a_1)<l(a),\ldots,l(a_k)<l(a).$ 
Then by inductive suggestion $a_1,\ldots,a_k\in\Gamma_1(\aleph).$ In other 
words, if $a_i=\alpha_{i,1}\ldots \alpha_{i,s_i},$ where 
$\alpha_{i,j}\in \aleph,$ then $\sum_{j=1}^{s_i}||\alpha_{i,j}||=1.$ Thus, 
$$a=\omega \alpha_{1,1}\ldots\alpha_{1,s_1}\ldots \alpha_{k,1}\ldots\alpha_{k,s_k},$$
and 
$$\mu_1(a)=||a||=||\omega||+\sum_{i=1}^k\sum_{j=1}^{s_i}||\alpha_{i,j}||=1-k +k =1.$$
By lemma \ref{1septe}, $\alpha_{i,s_i}\in {\frak X},$ for any $i=1,\ldots,k.$ 
Therefore, $\mu_i(a)\ge 1,$ for any $i\le \sum_{j=1}^ks_j.$ 
So, $a\in \Gamma_1(\aleph).$

Prove now $\Gamma_1)(\aleph)\subseteq \bar\Gamma_1(\aleph).$
By induction on $l(a)$ prove that any $a\in\Gamma_1(\aleph)$ can be 
presented in the form 
$a=x\in {\frak X},$ 
if $l(a)=1,$ or 
$a=\omega a_1\ldots a_r,$ where $||\omega||=1-r, a_1,\ldots,a_r\in 
\Gamma_1(\aleph),$ if $l(a)>1.$  

If $l(a)=1,$ then the statement is trivial. 
Suppose that our statement is true for elements of $\Gamma_1(\aleph)$ with 
length $<k$ and $a=\alpha_1\ldots\alpha_k, \alpha_i\in\aleph, i=1,2,\ldots,k.$ 
Let $\lambda_i=\mu_i(a).$ 

Suppose that $\mu_k=1\le \mu_{k-1}\le \cdots\le\mu_{l+1},$ but 
$\mu_{l}<\mu_{l+1}.$ This means that $\alpha_k,\ldots,\alpha_{l+1}\in {\frak X}$ and $\alpha_l\in\Omega.$ Let $|\alpha_l|=q>0.$ Then  
$q \le k-l,$ since $$\mu_l=1-q +\mathop{\underbrace{1+\cdots +1}}\limits_{k-l}\ge 1.$$  So, we can consider the element  
$c=\alpha_l\alpha_{l+1}\ldots \alpha_{l+q}\in \Gamma(\aleph).$ 
The word $c$ is a subword of $a.$ 
Moreover, $$||c||=||\alpha_l||+\ldots+||\alpha_{l+q}||=1-q+  
\mathop{\underbrace{1+\cdots +1}}\limits_{q}=1$$
and 
$$||\alpha_s\alpha_{s+1}\ldots\alpha_{l+q}||\ge 1,$$ 
for any $s=l+q,l+q-1,\ldots,l.$ 
So, $c\in\Gamma_1(\aleph).$ By inductive suggestion $c$ is an  $\Omega$-word.  
Therefore,  the word 
$$b=\beta_1\ldots\beta_{k-q-1}\in \Gamma(\aleph),$$
where
$$\beta_1=\alpha_1,\ldots,\beta_{l-1}=\alpha_{l-1},\beta_{l+1}=\alpha_{l+q+1},\ldots,
\beta_{k-q-1}=\alpha_k,$$
and $\beta_l\in {\frak X},$ has the following properties:  
$$l(b)=(l-1)+1+(k-l-q)=k-q-1< k,$$
$$||b||=||\alpha_1||+\cdots+||\alpha_{l-1}||+||\beta_l||+||\alpha_{l+q+1}||+\cdots+||\alpha_k||=||a||,$$
and
$$||\beta_s\ldots\beta_{k-q-1}||\ge 1,$$
for any $s=k-q-1,k-q-2,\ldots,1.$ These mean that $b\in\Gamma_1(\aleph)$ and 
$l(b)<l(a).$ By inductive suggestion, $b\in\bar\Gamma_1(\aleph),$ and 
$$b=\omega b_1\ldots b_r$$
for some $\omega\in\Omega$ and $b_1,\ldots,b_r\in\bar\Gamma_1(\aleph).$ 
Since $l(b_1),\ldots,l(b_r)<k,$ by inductive suggestion $b_1,\ldots,b_r\in\Gamma_1(\aleph).$ One of $b_s,$ where $1\le s\le r,$ contains $\beta_l.$ Instead of $\beta_l$ we can take the $\Omega$-word $c$ and obtain the $\Omega$-word that 
is equal to  $a.$ 
So, we established that $a\in\bar\Gamma_1(\aleph).$ 

Our theorem is proved. 

\begin{crl}\label{23may} For any $\Omega$-word $a,$  
$$\sum_{\omega\in \Omega\cap a}|\omega|=\omega deg\,a+xdeg\,a-1.$$
\end{crl}

{\bf Proof.} 
If $l(a)=1,$ then  $a\in {\frak X},$ and 
$\omega deg\,a=0,$ $xdeg\,a=1.$ So,
$\sum_{\omega\in \Omega\cap a}|\omega|=0=\omega deg\,a+xdeg\,a-1.$

Suppose that for $a,$ with $l(a)<k,$ the statement is true. Let 
$l(a)=k>1.$ 
By theorem \ref{22may}, any $\Omega$-word $a$ with length $k>1$ 
can be presented in the form $a=\eta a_1\ldots a_r,$ for some 
$\Omega$-words $a_1,\ldots,a_r$ with length $<k$ and some $\eta\in \Omega$ with 
$|\eta|=r.$   
Then 
$$\Omega\cap a= \{\eta\}\cup \cup_{i=1}^r\{\Omega\cap a_i\}.$$
Thus, 
$$\sum_{i=1}^r\omega deg\,a_i=\omega deg\,a-1.$$
By inductive suggestion,
$$
\sum_{\omega\in \Omega\cap a_i}|\omega|=\omega deg\,a_i+xdeg\,a_i-1,
$$
for $i=1,\ldots,r.$ 
Therefore, 
$$\sum_{\omega\in \Omega\cap a}|\omega|=$$
$$|\eta|+\sum_{i=1}^r\sum_{\omega\in \Omega\cap a_i}|\omega|=$$
$$r+\sum_{i=1}^r(\omega deg\,a_i+xdeg\,a_i-1)=$$
$$r+\sum_{i=1}^r\omega deg\,a_i+\sum_{i=1}^rxdeg\,a_i-r=$$
$$\omega deg\,a+xdeg\,a-1.$$

\begin{crl}\label{2sept} Let $\omega deg_i\,a$ 
be the number of entries of $\omega\in\Omega$ with $|\omega|=i$ in $a.$   
Then for any $\Omega$-word $a,$ 
$$xdeg\,a=\sum_{i\ge 1}(i-1)\omega deg_i\,a+1.$$
\end{crl}

{\bf Proof.} This is another formulation of corollary \ref{23may}. 

\begin{crl}\label{may23} Assume that  $\Omega$ consists of one element 
$\omega$ with $|\omega|=k.$  Then for any $\Omega$-word $a,$
$$xdeg\,a= 1+(k-1)\omega deg\,a.$$
\end{crl}

{\bf Proof.} It follows from  corollary \ref{2sept} and from 
the following facts:   
$\omega deg_k\,a=\omega deg\,a,$ and $\omega deg_i\,a=0,$ 
if $i\ne k.$ 

\section{$\Omega$-polynomials, $\Omega$-algebras and $\Omega$-identities}

{\bf Definition.} A linear combination of $\Omega$-words is called {\it
$\Omega$-polynomial.} Polynomial of the form $\lambda_a a,$ where
$a$ is a word and $\lambda_a\in K,$ is called  {\it a
monomial}. The monomial is called nontrivial, if $\lambda_a\ne 0.$ 
We say that $\lambda_a a$ is a part of $f$ or $\lambda_a a$ is monomial of 
$f,$ if $\lambda_a\ne 0.$  A space of $\Omega$-polynomials is denoted by
$K[{\Omega,\frak X}].$

Let $U$ and $M$ be some  vector spaces. 
Denote by $T^k(U,M)$ the space of polylinear maps 
$\psi:\mathop{\underbrace{U\times \cdots \times U}}\limits_{k}
\rightarrow M,$ if $k>0,$ $T^0(U,M)=M,$ and 
$T^k(U,M)=0,$ if $k<0.$ Let $T^*(U,M)=\oplus_kT^k(U,M).$ 
If $\psi\in T^k(U,M),$ we will write $|\psi|=k.$ 

Let $\wedge^k(U,M)$ be the subspace of $T^k(U,M)$ consisting of 
skew-symmetric maps, $\wedge^0(U,M)=M,$ $ \wedge^k(U,M)=0,$ if $k<0$  
and $\wedge^*(U,M)=\oplus_k\wedge^k~(U,M).$ 

Let $\Omega=\{\omega_1,\omega_2,\ldots\}$ be some alphabet with an  
arity map $|\;,\;|: \Omega\rightarrow {\bf Z}^+.$
Suppose that  to each $\omega\in\Omega$ one corresponds some  
homogeneous map $\omega_U\in T^{|\omega|}(U,U).$ 
In this case we will say that $U$  has a structure of 
$\Omega\mbox{-\it algebra}$.

Notice that for any $\Omega$-algebra $U$ and for any $\Omega$-word $a$ 
one can make substitutions for its parameters $x_i\mapsto u_i\in U$ and $\omega_i\mapsto \omega_{iU}.$ The easy way to see it is by presenting $a$ in the form
$a=\omega a_1\ldots a_k,$ where $|\omega|=k$ and 
$a_1,\ldots,a_k$ are words of smaller degree than the $deg\,a.$ 
One can assume that in $a_1,\ldots,a_k$ our substitutions are correctly 
defined. Then $a$ would be correctly defined also. 

So, for any polynomial $f\in K[\Omega,\frak X]$ we can make substitutions in 
its parameters by elements of $U$ and operations on $U.$ 
If $f$ depends, say, from
parameters $x_1,\ldots,x_k, \omega_1,\ldots,\omega_l,$ then we obtain a  
correctly defined element 
$f_U=f(u_1,\ldots,u_k,\omega_{1U},\ldots,\omega_{lU})\in  U$ 
for any $u_1,\ldots,u_k\in U.$

{\bf Definition.} The polynomial $f\in K[{\Omega,\frak X}]$ is called an 
{\it $\Omega$-polynomial identity}, or simply, $\Omega$-identity on $\Omega$-algebra $U,$ if $$f_U(u_1,\ldots,u_k,\omega_{1U},\ldots,\omega_{lU})=0,$$ for 
any $u_1,\ldots,u_k\in U.$ 

Let us given two $\Omega$-polynomial identities $f=0$ and $g=0.$ 
We say that the identity $f=0$ {\it follows} from the identity $g=0,$ 
and denote $g=0\Rightarrow f=0,$
if $f_U=0$ for any $\Omega$-algebra $U,$ 
such that $g_U=0.$ The identities $f=0$ and $g=0$ are called 
{\it $\Omega$-equivalent,} or just equivalent, if $f=0\Rightarrow g=0$ and $g=0\Rightarrow f=0.$ 

Further, to simplify denotions we will often identify 
the polynomial $f=f(t_1,\ldots,t_k,\omega_1,\ldots,\omega_l)$ by the 
result of substitution $f_U=f(u_1,\ldots,u_k,\omega_{1U},\ldots,\omega_{lU})
\in  U$ and call $f(u_1,\ldots,u_k,\omega_{1U},\ldots,\omega_{lU})$ 
shortly as a $\Omega$-polynomail, or just a polynomial.
Notice that a formal definiton of  $\Omega$-words does not need  
any brackets and comma's, but for practical use it is more convenient 
to use brackets and comma's. We will 
use brackets keeping in mind that we will do it from the right to 
the left as in the proof of theorem \ref{22may}.

{\bf Example.} Let $\Omega=\{\omega_3,\omega_2,\omega_2':
|\omega_3|=3, |\omega_2|=|\omega_2'|=2.$ 
Let 
$a=\omega_3\omega_2x_1x_2x_3\omega_2'x_4x_5\in \Gamma_1(\aleph),$
Then for any $\Omega$-algebra $U$ 
and  for any $u_1,\ldots,u_5\in U,$ 
$$a_U=\omega_{3U}(\omega_{2U}(u_1,u_2),u_3,\omega_{2U}'(u_4,u_5))\in U,$$
or simply,
$$a=\omega_{3}(\omega_{2}(x_1,x_2),x_3,\omega_{2}'(x_4,x_5)).$$

\section{$3$-Lie  algebras}

\begin{thm} \label{27august} If $p=0$ or  $p>3,$ then any $3$-Lie-Poisson algebra is strong.
\end{thm}

{\bf Proof.} Let $(U,\cdot,\omega)$ be $3$-Lie-Poisson. 
Recall that 
$$r^{\omega}(u_1,\ldots,u_6)=$$
$$\omega(u_1,u_2,u_3)\cdot\omega(u_4,u_5,u_6-
\omega(u_1,u_2,u_5)\cdot \omega(u_4,u_3,u_6)$$
$$+\omega(u_1,u_2,u_6)\cdot\omega(u_4,u_3,u_5)+
\omega(u_4,u_2,u_3)\cdot\omega(u_1,u_5,u_6)$$ $$
-\omega(u_4,u_2,u_5)\cdot \omega(u_1,u_3,u_6)+\omega(u_4,u_2,u_6)\cdot\omega(u_1,u_3,u_5)$$
is symmetric in two arguments $u_1$ and $u_4$ 
and  skew-symmetric in three arguments $u_3,u_5,u_6.$  

By theorem 6 \cite{DalTakh}, for $3$-Lie-Poisson algebras, $r^{\omega}(u_1,\ldots,u_6)=0,$ for any $u_1,\ldots,u_6\in U.$ 
One can check that  
$$3f^{\omega}(u_1,u_2,u_3,u_4,u_5,u_6)=$$
$$2 r^{\omega}(u_1,u_2,u_3,u_4,u_5,u_6)+r^{\omega}(u_2,u_3,u_1,u_4,u_5,u_6)$$
$$-r^{\omega}(u_2,u_4,u_1,u_3,u_5,u_6)+r^{\omega}(u_2,u_5,u_1,u_3,u_4,u_6)-r^{\omega}(u_2,u_6,u_1,u_3,u_4,u_5).$$
So, $f^{\omega}=0$ is also identity for $(U,\cdot,\omega),$ if  $p\ne 3.$

\begin{prp}\label{22sept} $f^{\omega}=0\Rightarrow r^{\omega}=0.$
\end{prp}

{\bf Proof.} One can check that:  
$$-r^{\omega}(u_1,u_n,u_2,\ldots,\hat{u_n},\ldots, u_{2n})= $$
$$f^{\omega}(u_1,u_2,\ldots,u_n,\ldots,u_{2n})+(-1)^n
f^{\omega}(u_2,\ldots,u_{n},u_1,u_n,\ldots,u_{2n}).$$
Therefore, the identity $r^{\omega}=0$ follows from the identity 
$f^{\omega}=0.$   

{\bf Remark.} 
Notice that for $n=2$ the identities $r^{\omega}=0$ and $f^{\omega}=0$ 
are different. More exactly, the identity $r^{\omega}=0$ does not appear, if 
$n=2.$ There exist 2-Lie-Poisson algebras that are not strong. 
Let us give an example of such algebras.  Let $K_{2l}=K[x_1,\ldots,x_{2l}]$ 
and $\omega=\sum_{i=1}^l\der_i\wedge\der_{i+l}.$ Then $(K_{2l},\cdot,\omega)$ 
is $2$-Lie-Poisson. 
It is easy to check that $(K_2,\cdot,\der_1\wedge\der_2)$ satisfies the identity $f^{\der_1\wedge\der_2}=0.$ If $l>1,$ the polynomial  
$$f^{\omega}(a,u,v,w)=\omega(a,u)\cdot\omega(v,w)+\omega(a,v)\cdot \omega(w,u)+
\omega(a,w)\cdot\omega(u,v)$$
is not an identity. For instance, 
$$f^{\omega}(x_1,x_2,x_3,x_4)=\omega(x_1,x_2)\cdot\omega(x_3,x_4)=1\ne 0.$$
So, the algebra $(K_{2l},\cdot,\omega)$ is strong $2$-Lie-Poisson, if and only if $l=1.$

\section{ Minimal identities for $3$-Jacobians}

\begin{thm}\label{minimal} 
$(p\ne 2, 3)$ Any polynomial identity of $\Omega$-degree $2$
of Jacobian algebra $(K[x_1,x_2,x_3],Jac_3^S)$ follows from $3$-Lie 
and skew-symmetric identities for $3$-multiplication $Jac_3^S.$ 

\end{thm}

{\bf Proof.} Let $\omega=Jac_3^S.$
Define polynomials $g^{\omega},$  $h^{\omega}$ and $q^{\omega}$ by  
$$g^{\omega}=g^{\omega}(t_1,\ldots,t_5)=$$
$$\omega t_1t_2\omega t_3t_4t_5-
\omega\omega t_1,t_2,t_3t_4t_5+$$
$$
\omega\omega t_1t_2t_4t_3t_5-
\omega\omega t_1t_2t_5t_3t_4,$$

$$h^{\omega}=h^{\omega}(t_1,\ldots,t_5)=$$
$$\omega t_1t_2\omega t_3t_4t_5-
\omega t_1t_3\omega t_2t_4t_5+$$
$$
\omega t_1t_4\omega t_2t_3t_5-
\omega t_1t_5 \omega t_2t_3t_4,$$

$$q^{\omega}=q^{\omega}(t_1,\ldots,t_5)=$$
$$\sum_{i<j<5}(-)^{i+j}\omega t_i t_j\omega t_1\ldots\hat{t_i}\ldots\hat{t_j}\ldots t_5.$$

Notice that 
$$3 h^{\omega}(t_1,t_2,t_3,t_4,t_5)=$$
$$\sum_{2\le i<j\le 5}(-1)^{i+j}g^{\omega}(t_i,t_j,t_1,\ldots,\hat{t_i},\ldots\,\hat{t_j},\ldots, t_5).$$
Further,
$$\sum_{1\le i<j\le 5}(-1)^{i+j}g^{\omega}(t_1,\ldots,t_5)=$$
$$2 \sum_{1\le i<j\le 5}(-1)^{i+j}\omega t_i t_j\omega t_1\ldots\hat{t_i}\ldots\hat{t_j}\ldots t_5.$$
Therefore,
$$2 q^{\omega}(t_2,t_3,t_4,t_5,t_1)=\sum_{1\le i<j\le 5}(-1)^{i+j}g^{\omega}(t_1,\ldots,t_5)-2 h^{\omega}(t_1,\ldots,t_5).$$
So, $h^{\omega}=0$ and  $q^{\omega}=0$ are identities on $(U,\omega)$, 
if $g^{\omega}=0$ is an identity and $p\ne 2,3.$ 

Let $f$ be any polynomial of $wdeg\,f= 2.$  Since $|\omega|=3,$ 
by corollary \ref{2sept} we should prove that any polynomial of the form 
$f(t_1,\ldots,t_5)=\sum_{i_1<i_2,i_3<i_4<i_5}
\lambda_{i_1i_2}\omega(t_{i_1},t_{i_2},\omega(t_{i_3},t_{i_4},t_{i_5})),$ 
such that 
$$f(u_1,\ldots,u_5)=0,$$
for any $u_1,\ldots,u_5\in U=K[x_1,x_2,x_3],$   
is a linear combination of polynomials that can be obtained from the 
polynomial $g^{\omega}$ by permutation of variables $t_1,\ldots,t_5.$ 

We have 
$$f(x_1,x_2,x_1,x_2,x_3^2)=0\Rightarrow \lambda_{12}-\lambda_{14}+\lambda_{34}-\lambda_{23}=0,$$
$$f(x_2,x_1,x_1,x_2,x_3^2)=0\Rightarrow \lambda_{12}+\lambda_{13}+\lambda_{24}+\lambda_{34}=0,$$
$$f(x_1,x_2,x_1,x_3^2,x_2)=0\Rightarrow \lambda_{12}+\lambda_{15}-\lambda_{23}-\lambda_{35}=0,$$
$$f(x_2,x_1,x_1,x_3^2,x_2)=0\Rightarrow \lambda_{12}+\lambda_{13}-\lambda_{25}-\lambda_{35}=0,$$
$$f(x_1,x_2,x_3^2,x_1,x_2)=0\Rightarrow \lambda_{12}+\lambda_{15}+\lambda_{24}
+\lambda_{45}=0,$$
$$f(x_2,x_1,x_3^2,x_1,x_2)=0\Rightarrow -\lambda_{12}+\lambda_{14}+\lambda_{25}-\lambda_{45}=0.$$
The obtained system of linear equations has rank $5$ and parameters 
$\lambda_{15},\lambda_{25},\lambda_{34},\lambda_{35}, \lambda_{45}
$ can be chosen as a free. So, any polynomial $f,$ 
such that $xdeg\,f=5$ and $f=0$ is an identity on $(K[x_1,x_2,x_3],\omega),$ is a linear combination of the following five polynomials 
$$f_{45}=-\omega t_1t_2\omega t_3t_4t_5+\omega t_1t_3\omega t_2t_4t_5-\omega t_2t_3 \omega t_1 t_4 t_5+\omega t_4 t_5 \omega t_1 t_2 t_3,$$
$$f_{35}=\omega t_1t_2\omega t_3t_4t_5+\omega t_1t_4\omega t_2t_3t_5-\omega t_2t_4 \omega t_1 t_3 t_5+\omega t_3 t_5 \omega t_1 t_2 t_4,$$
$$f_{34}=\omega t_1t_2\omega t_3t_4t_5-\omega t_1t_3\omega t_2t_4t_5+\omega t_1t_4 \omega t_2 t_3 t_5$$
$$+\omega t_2 t_3  \omega t_1 t_4 t_5-\omega t_2 t_4\omega t_1 t_3 t_5+\omega t_3 t_4 \omega t_1 t_2 t_5,$$
$$f_{25}=\omega t_1t_2\omega t_3t_4t_5+\omega t_2t_3\omega t_1t_4t_5-
\omega t_2t_4 \omega t_1 t_3 t_5+\omega t_2 t_5 \omega t_1 t_3 t_4,$$
$$f_{15}=-\omega t_1t_2\omega t_3t_4t_5+\omega t_1t_3\omega t_2t_4t_5-\omega t_1t_4 \omega t_2 t_3 t_5+\omega t_1 t_5 \omega t_2 t_3 t_4.$$
We see that 
$$f_{45}=g^{\omega}(t_4,t_5,t_1,t_2,t_3),$$
$$f_{35}=g^{\omega}(t_3,t_5,t_1,t_2,t_4),$$
$$f_{34}=-q^{\omega}, $$
$$f_{25}=-h^{\omega}(t_2,t_1,t_3,t_4,t_5),$$
$$f_{15}=-h^{\omega}.$$
So, any identity of $\frak X$-degree 5 follows from $3$-Lie 
identity $g^{\omega}=0$ and skew-symmetric identity. 

{\bf Conjecture.} Let $p=0$ and $n>2.$ 
Any identity of $\Omega$-degree $\le 3$ of $n$-Lie 
algebra $(K[x_1,\ldots,x_n],Jac_n^S)$ follows 
from $n$-Lie and skew-symmetric identities.   
Any identity of $\Omega$-degree $\le 3$ of $n$-Lie-Poisson algebra 
$(K[x_1,\ldots,x_n],\cdot, Jac_n^S)$ follows 
from $n$-Lie and skew-symmetric identities for Jacobian,
 Leibniz rule for $\cdot$ and the identity $f^{Jac_n^S}=0.$ 

Here we suppose that commutativity and associativity identities for 
a binary operation  are given. Definitions of $n$-Lie-Poisson algebras and the 
polynomial $f^{\omega}$ are given in the next section.  

\section{Constructions of $n$-Lie algebras by $n$-Lie-Poisson algebras}

Usually $n$-Lie algebras are considered for $n>1.$ 
Complete this definition  by considering the case  $n=1.$ 
Call any vector space $U$ with a linear map  $f:U\rightarrow U,$ i.e.,  
$(U,f)$ a $1$-Lie algebra.

{\bf Definition}. Let $A=(U,\cdot,\omega)$ be an algebra with two operations: 
$(U,U)\rightarrow U, (u,v)\mapsto u\cdot v$ be a bilinear multiplication and 
$\omega: \wedge ^n U\rightarrow U$ be a skew-symmetric $n$-linear 
multiplication. We say that $A$ is {\it $n$-Lie-Poisson,} if 
\begin{itemize}
\item $(U,\cdot)$ is an associative commutative algebra
\item $(U,\omega)$ is an $n$-Lie algebra.
\item $\omega(u\cdot u', u_2,\ldots,u_n)=\omega(u,u_2,\ldots,u_n)\cdot u'+
u\cdot \omega(u',u_2,\ldots,u_n),$ for any $u,u',u_2,\ldots,u_n\in U.$ 
\end{itemize}

Consider on $U$ one more polynomial 
$$f^\omega(u_1,\ldots,u_{n-1},v_1,\ldots,v_n)=$$
$$\sum_{i=1}^{n+1}(-1)^i\omega(u_1,\ldots,u_{n-1},v_{i})\cdot 
\omega(v_{1},\ldots,\hat{v_i},\ldots,v_{n+1}).
$$

Let us given an algebra $(U,\cdot,\omega)$ with 
one binary operation $(u,v\mapsto u\cdot v$ and one $n$-ary operation 
$\wedge^nU\rightarrow U, (u_1,\ldots,u_n)\mapsto \omega(u_1,\ldots,u_n).$  

{\bf Definition.} An $n$-Lie-Poisson algebra is called {\it strong}, if it satisfies the identity $f^{\omega}=0.$ 

{\bf Example.} We will see below that the Jacobian algebra $(U,\cdot,Jac_n^S)$ 
is strong $n$-Lie-Poisson and that 
the algebra $(U,\cdot,Jac_{n+1}^W)$ 
satisfies the identity $f^{\omega}=0.$ 

Let $(U,\cdot,\omega)$ be $n$-Lie-Poisson algebra.  Recall that 
$D\in Der\,(U,\cdot,\omega),$ if  
\begin{itemize}
\item  $D\in Der\,(U,\cdot),$ i.e., $D(u\cdot v)=D(u)\cdot v+u\cdot D(v),$ 
\item $D\in Der\,(U,\omega),$ i.e., 
$$D(\omega(u_1,\ldots,u_n))=\sum_{i=1}^n\omega(u_1,\ldots,u_{i-1},
D{u_i},u_{i+1},\ldots,u_n).$$
\end{itemize}
Call in a such case $D$ a $n$-Lie-Poisson derivation.

{\bf Example.} Let $n=1.$ Then $(U,\cdot,\omega)$ is $n$-Lie-Poisson, if 
$\omega: U\rightarrow U$ is a derivation of $(U,\cdot).$ 
If $D$ is an $n$-Lie-Poisson derivation of $(U,\cdot,\omega),$ then 
$[D,\omega]=0.$ 

\begin{thm} \label{14august1} Let $(U,\cdot,\omega)$ be a strong 
$n$-Lie-Poisson algebra and $D\in Der\,(U,\cdot,\omega).$ 
Construct on $U$ a new skew-symmetric 
$(n+1)$-multiplication $\bar \omega=D\wedge \omega.$ 
Then $(U,\cdot,\bar\omega)$ is  strong $(n+1)$-Lie-Poisson. 
\end{thm}

{\bf Proof.} We have 
$$\bar\omega(u\cdot u',u_1,\ldots,u_{n})=$$
$$ D(u\cdot u')\cdot\omega (u_1,\ldots,u_n)+\sum_{i=1}^n(-1)^iD(u_i)\cdot\omega(u\cdot u',u_1,\ldots,\hat{u_i},\ldots,u_n)= $$

$$ (D(u)\cdot u')\cdot\omega (u_1,\ldots,u_n)+
u\cdot D(u')\cdot\omega (u_1,\ldots,u_n)
$$
$$+\sum_{i=1}^n(-1)^iD(u_i)\cdot u\cdot \omega(u',u_1,\ldots,\hat{u_i},\ldots,u_n)
+$$
$$\sum_{i=1}^n(-1)^iD(u_i)\cdot \omega(u,u_1,\ldots,\hat{u_i},\ldots,u_n)\cdot u'= $$

$$ (D(u)\cdot\omega (u_1,\ldots,u_n))\cdot u'
+\sum_{i=1}^n(-1)^i(D(u_i)\cdot \omega(u,u_1,\ldots,\hat{u_i},\ldots,u_n))\cdot u'+$$
$$u\cdot (D(u')\cdot\omega (u_1,\ldots,u_n))
+\sum_{i=1}^n(-1)^iu\cdot(D(u_i)\cdot\omega(u',u_1,\ldots,\hat{u_i},\ldots,u_n))=$$

$$D\wedge \omega(u,u_1,\ldots,u_n)\cdot u'+u\cdot \omega(u',u_1,\ldots,u_n)$$

Further,
$$\sum_{i=1}^{n+2}(-1)^iD\wedge \omega(u_1,\ldots,u_n,v_i)\cdot D\wedge \omega(v_1,\ldots,\hat{v_i},\ldots,v_{n+2})=$$
$$X_1+X_2,$$
where 
$$X_1=\sum_{i=1}^{n+2}\sum_{j=1}^n\sum_{1\le s\le n+2,s\ne i}
(-1)^{i+j+s+\delta(s>i)}\cdot$$
$$
D(u_j)\cdot D(v_s)\cdot \omega(u_1,\ldots,\hat{u_j},\ldots, u_n,v_i)
\cdot \omega(v_1,\ldots,\hat{v_s},\ldots,\hat{v_i},\ldots,v_{n+2}),$$

$$X_2=
\sum_{i=1}^{n+2}\sum_{1\le s\le n+2, s\ne i}
(-1)^{i+n+s+\delta(s>i)+1}\cdot 
$$
$$D(v_i)\cdot D(v_s)\cdot
 \omega(u_1,\ldots,u_n)\cdot \omega(v_1,\ldots,\hat{v_s},\ldots,\hat{v_i},\ldots,v_{n+2}).$$
Notice that, by identity $f^{\omega}=0,$ 
$$X_1=\sum_{j=1}^n\,\sum_{1\le s\le n+2}
(-1)^{j+s}D(u_j)\cdot D(v_s)\cdot $$
$$
\left(\sum_{1\le i\le n+2, i\ne s}(-1)^{i+\delta(s>i)}\omega(u_1,\ldots,\hat{u_j},\ldots, u_n,v_i)
\cdot \omega(v_1,\ldots,\hat{v_s},\ldots,\hat{v_i},\ldots,v_{n+2})\right)$$
$$
=0$$
and by commutativity of multiplication $\cdot,$ 
$$X_2=
(-1)^{n}\omega(u_1,\ldots,u_n)\cdot
$$
$$
\left(\sum_{i=1}^{n+2}\sum_{i< s\le n+2}
(-1)^{i+s}D(v_i)\cdot D(v_s)\cdot 
\omega(v_1,\ldots,\hat{v_i},\ldots,\hat{v_s},\ldots,v_{n+2})\right.
-$$
$$
\left.\sum_{i=1}^{n+2}\sum_{i< s\le n+2}
(-1)^{i+s}D(v_s)\cdot D(v_i)\cdot 
\omega(v_1,\ldots,\hat{v_i},\ldots,\hat{v_s},\ldots,v_{n+2})
\right)=0$$
So, $f^{\bar\omega}=0$ is the identity. 

Notice that  
$$D\wedge\omega(u_1,\ldots,u_n,D\wedge \omega(u_{n+1},\ldots,u_{2n+1}))$$
$$-\sum_{i=1}^{n+1}(-1)^{i+n+1}D\wedge\omega(D\wedge 
\omega(u_1,\ldots,u_n,u_{n+i}),u_{n+1},\ldots,\widehat{u_{n+i}},\ldots,u_{2n+1})
=$$

$$ \sum_{i=1}^{2n+1}D(D(u_i))\cdot g_i(u_1,\ldots,\hat{u_i},\ldots,u_{2n+1})+$$
$$\sum_{1\le i,j\le 2n+1, i\ne j}D(u_i)\cdot D(u_j)\cdot 
g_{i,j}(u_1,\ldots,\hat{u_i},\ldots,\hat{u_j},\ldots,u_{2n+1})+$$
$$\sum_{i=1}^{2n+1}\sum_{j=1}^{2n+1}D(u_i)\cdot h_{i,j}(u_1,\ldots,\hat{u_i},\ldots,u_{j-1},D(u_j),u_{j+1},\ldots,u_{2n+1}),$$
for some polynomials $g_i, g_{i,j}$ and $h_{i,j},$   
that do not depend from $D.$ 

The following relations can be obtained by tedious calculations.   

If $n<i\le 2n+1,$ then  $g_i=0.$
If  $1\le i\le n,$ then by the identity $f^{\omega}=0,$ 
$$g_i(u_1,\ldots,\hat{u_i},\ldots,u_{2n+1})=$$
$$(-1)^{i+n+1}\sum_{j=n+1}^{2n+1}
(-1)^{j}\omega(u_1,\ldots,\hat{u_i},\ldots,u_n,u_j)\cdot 
\omega(u_{n+1},\ldots,\hat{u_j},\ldots,u_{2n+1})
=$$
$$(-1)^{i+n+1}f^{\omega}(u_1,\ldots,\hat{u_i},\ldots,u_{2n+1})=0.$$

If $1\le i,j\le n,$ or $n<i,j\le 2n+1,$ then $g_{i,j}=0.$
If $1\le i\le n,$ $n<j\le  2n+1,$ by $n$-Lie identity $g_{i,j}=0.$

If $n<j\le 2n+1,$ then $h_{i,j}=0.$ If $1\le j\le n,$ then 
$$h_{i,j}=\pm 
f^{\omega}(u_1,\ldots,\hat{u_i},\ldots,u_{j-1},D(u_j),u_{j+1},\ldots,u_n,
u_{n+1},\ldots,u_{2n+1}),$$
in the case of $1\le i\le n,$ and 
$$h_{i,j}=\pm 
f^{\omega}(u_1,\ldots,\hat{u_{j}},\ldots,u_n,
u_{n+1},\ldots,u_{i-1},D(u_i),u_{i+1},\ldots,u_{2n+1}),$$
in the case of $n<i\le 2n+1.$ Thus, in both cases, by the identity $f^{\omega}=0,$ we have $h_{i,j}=0.$ 

So, the multiplication $\bar\omega$ satisfies $(n+1)$-Lie identity.   

\begin{crl} \label{20aug} 
Let $U$ be an associative commutative algebra with commuting derivations $\der_1,\ldots,\der_n.$ Let $\omega=\der_1\wedge \cdots\wedge \der_n$ 
Then $(U,\cdot,\omega)$ is strong $n$-Lie-Poisson. 
\end{crl}

The algebra $(U,\cdot,\omega)$ constructed in corollary~\ref{20aug}  
is called  {\it $n$-Lie Jacobian algebra of type $S.$}

{\bf Proof.} We will argue by induction on $n$. If $n=1,$ then $(U,\cdot,\der_1)$ is $1$-Lie-Poisson, if and only if $\der_1\in Der\,(U,\cdot).$  
In this case the identity $f^{\omega}=0$ follows from the commutativity 
law for the algebra $(U,\cdot).$ 

Suppose that the algebra $(U,\cdot,\der_1\wedge\cdots\wedge\der_{n-1})$ is 
strong $(n-1)$-Lie-Poisson. By definition $\der_n\in Der\,(U,\cdot).$ Since 
$[\der_n,\der_i]=0,$ for any $i=1,\ldots,n-1,$ then 
$$\der_n(\der_1\wedge\cdots\wedge\der_{n-1}(u_1,\ldots,u_{n-1})=$$
$$\der_n(\sum_{\sigma\in Sym_{n-1}}sign\,\sigma\,\der_{\sigma(1)}u_1\cdot \cdots \cdot \der_{\sigma(n-1)}u_{n-1})=$$
$$\sum_{i=1}^{n-1}\sum_{\sigma\in Sym_{n-1}}sign\,\sigma\,\der_{\sigma(1)}u_1\cdot \cdots\cdot
\der_{\sigma(i-1)}u_{i-1}\der_n(\der_{\sigma(i)}u_i)\cdot\der_{\sigma(i+1)}u_{i+1}\cdot \cdots  \cdot \der_{\sigma(n-1)}u_{n-1}=
$$
$$\sum_{i=1}^{n-1}\sum_{\sigma\in Sym_{n-1}}sign\,\sigma\,\der_{\sigma(1)}u_1\cdot \cdots\cdot
\der_{\sigma(i-1)}u_{i-1}\cdot\der_{\sigma(i)}\der_n(u_i)\cdot\der_{\sigma(i+1)}u_{i+1}\cdot \cdots  \cdot \der_{\sigma(n-1)}u_{n-1}=
$$
$$\sum_{i=1}^{n-1}\der_1\wedge\cdots\wedge\der_{n-1}(u_1,\ldots,u_{i-1},\der_nu_i,u_{i+1},\ldots,u_{n-1}).$$
So, $\der_n\in Der\,(U,\der_1\wedge\cdots\wedge\der_{n-1}).$ 
Thus we can apply  theorem \ref{14august1}. By this theorem we obtain that 
$(U,\cdot,\der_1\wedge\cdots\wedge\der_n)$ is strong $n$-Lie-Poisson.

\begin{thm}\label{14august2} Let $(U,\cdot,\omega)$ be strong 
$n$-Lie-Poisson. Endow $U$ by a new skew-symmetric $(n+1)$-multiplication 
$\tilde\omega=
id\wedge \omega,$ where $id:U\rightarrow U, u\mapsto u,$ is the identity 
map. Then $(U,\tilde\omega)$ is $(n+1)$-Lie algebra. It satisfies the identity 
$f^{\tilde\omega}=0.$ If $U$ has unit $1,$ then it satisfies one more 
identity 
\begin{equation}\label{19aug}
\tilde\omega(u\cdot u',u_1,\ldots,u_n)-u\cdot\tilde \omega(u',u_1,\ldots,u_n)
-u'\cdot \tilde\omega(u,u_1,\ldots,u_n)=
\end{equation}
$$-(u\cdot u')\cdot i(1)\tilde\omega(u_1,\ldots,u_n)
$$
\end{thm}

{\bf Proof.} Since $L_{u_2,\ldots,u_n}\in Der\,(U,\cdot),$ and $1\cdot1=1,$
$$i(1)\omega=0.$$
 Therefore, 
\begin{equation}\label{19augu}
i(1)(id\wedge\omega)=\omega
\end{equation} 
So, 
for any $u,u',u_1,\ldots,u_n\in U,$ 
$$\tilde\omega(u\cdot u',u_1,\ldots,u_n)=Y_1+Y_2,$$
$$u\cdot\tilde \omega(u',u_1,\ldots,u_n)=Z_1+Z_2,$$
$$u'\cdot \tilde\omega(u,u_1,\ldots,u_n)=W_1+W_2,$$
where
$$Y_1=(u\cdot u')\cdot\omega(u_1,\ldots,u_n),$$
$$Y_2=\sum_{i=1}^n(-1)^iu_i\omega(u\cdot u',u_1,\ldots,\hat{u_i},\ldots,u_n),$$
$$Z_1=u\cdot (u'\cdot \omega(u_1,\ldots,u_n)),$$
$$Z_2=\sum_{i=1}^n(-)^iu\cdot 
(u_i\cdot\omega(u',u_1,\ldots,\hat{u_i},\ldots,u_n)),$$
$$W_1=u'\cdot (u\cdot \omega(u_1,\ldots,u_n)),$$
$$W_2=\sum_{i=1}^n(-)^iu'\cdot 
(u_i\cdot\omega(u,u_1,\ldots,\hat{u_i},\ldots,u_n)).$$
Since $L_{u_1,\ldots,\hat{u_i},\ldots,u_n}\in Der\,(U,\cdot),$ 
$$Y_1=Z_1=W_1, \quad Y_2=Z_2+W_2.$$
Therefore, by (\ref{19augu}), the identity (\ref{19aug}) is true. 

We have 
$$f^{\tilde\omega}(u_1,\ldots,u_n,v_1,\ldots,v_{n+2})=$$
$$\sum_{i=1}^{n+2}(-1)^iid\wedge\omega(u_1,\ldots,u_n,v_i)\cdot \tilde\omega(v_1,\ldots,\hat{v_i},\ldots,v_{n+2})=$$

$$\sum_{i=1}^{n+2}\sum_{j=1}^n(-1)^{i+j}
u_j\cdot\omega(u_1,\ldots,\hat{u_j},\ldots,u_n,v_i)
\cdot \tilde\omega(v_1,\ldots,\hat{v_i},\ldots,v_{n+2})+$$
$$\sum_{i=1}^{n+2}(-1)^{i+n+1}
v_i\cdot\omega(u_1,\ldots,u_n)
\cdot \tilde\omega(v_1,\ldots,\hat{v_i},\ldots,v_{n+2})=$$
$$T_1+T_2,$$
where
$$T_1=\sum_{i=1}^{n+2}\sum_{j=1}^n\sum_{1\le s\le n+2,s\ne i}(-1)^{i+j+s+\delta(s>i)}\cdot$$
$$
u_j \cdot v_s\cdot\omega(u_1,\ldots,\hat{u_j},\ldots,u_n,v_i)
\cdot \omega(v_1,\ldots,\hat{v_s},\ldots,\hat{v_i},\ldots,v_{n+2}),$$

$$T_2=\sum_{i=1}^{n+2}\sum_{1\le s\le n+2, s\ne i}
$$
$$(-1)^{i+s+\delta(s>i)+n+1}
v_i\cdot v_s\cdot\omega(u_1,\ldots,u_n)
\cdot \omega(v_1,\ldots,\hat{v_s},\ldots,\hat{v_i},\ldots,v_{n+2}).$$
Since $(U,\cdot)$ is commutative,
$$T_2=
(-1)^{n+1}\cdot\omega(u_1,\ldots,u_n)\cdot$$
$$ \left(\sum_{i=1}^{n+2}\sum_{1\le s\le n+2, s\ne i}(-1)^{i+s+\delta(s>i)}
v_i\cdot v_s\cdot \omega(v_1,\ldots,\hat{v_s},\ldots,\hat{v_i},\ldots,v_{n+2})\right)=0.$$

\medskip

\noindent Notice that 
$$T_1=\sum_{j=1}^n\sum_{s=1}^{n+2}(-1)^{j+s}u_j \cdot v_s\cdot$$
$$
\left(
\sum_{1\le i\le n+2, i\ne s}(-1)^{i+\delta(s>i)}
\omega(u_1,\ldots,\hat{u_j},\ldots,u_n,v_i)
\cdot \omega(v_1,\ldots,\hat{v_s},\ldots,\hat{v_i},\ldots,v_{n+2})\right)=$$
$$\sum_{j=1}^n\sum_{s=1}^{n+2}(-1)^{j+s}u_j \cdot v_s\cdot f^{\omega}(u_1,\ldots,\hat{u_j},\ldots,u_n,v_1,\ldots,\hat{v_s},\ldots,v_{n+2}).$$
Since $f^{\omega}=0$ is identity on $(U,\cdot,\omega),$ we have $T_1=0.$ 
So, the algebra $(U,\cdot,\tilde\omega)$ satisfies the identity 
$f^{\tilde\omega}=0.$

We have
$$id\wedge\omega(u_1,\ldots,u_n,id\wedge \omega(u_{n+1},\ldots,u_{2n+1}))
=$$

$$\sum_{i=n+1}^{2n+1}(-1)^{i+n+1}id\wedge \omega(u_1,\ldots,u_n,
u_{i}\cdot \omega(u_{n+1},\ldots,\hat{u_i},\ldots,u_{2n+1}))=$$

$$\sum_{j=1}^n\sum_{i=n+1}^{2n+1}(-1)^{i+j+n}u_j\cdot  
\omega(u_1,\ldots,\hat{u_j},\ldots,u_n,
u_{i}\cdot \omega(u_n,\ldots,\hat{u_i},\ldots,u_{2n+1}))+$$
$$\sum_{i=n+1}^{2n+1}(-1)^{i+1}u_i\cdot \omega(u_1,\ldots,u_n)\cdot 
\omega(u_{n+1},\ldots,\hat{u_i},\ldots,u_{2n+1})
=$$
$$A_1+A_2+A_3,$$ 
where
$$A_1=\sum_{j=1}^n\sum_{i=n+1}^{2n+1}(-1)^{i+j+n}u_j\cdot  
\omega(u_1,\ldots,\hat{u_j},\ldots,u_n,
u_{i})\cdot \omega(u_n,\ldots,\hat{u_i},\ldots,u_{2n+1})),$$
$$A_2=\sum_{j=1}^n\sum_{i=n+1}^{2n+1}(-1)^{i+j+n}u_ju_i\cdot  
\omega(u_1,\ldots,\hat{u_j},\ldots,u_n,
\omega(u_n,\ldots,\hat{u_i},\ldots,u_{2n+1})),$$
$$A_3=\omega(u_1,\ldots,u_n)\cdot 
\sum_{i=n+1}^{2n+1}(-1)^{i+1}u_i\cdot  
\omega(u_{n+1},\ldots,\hat{u_i},\ldots,u_{2n+1}). $$

On the other hand,
$$\sum_{i=n+1}^{2n+1}
(-1)^{i+n+1}(id\wedge\omega)(id\wedge \omega(u_1,\ldots,u_n,u_i),
u_{n+1},\ldots,\hat{u_i},\ldots,u_{2n+1})=$$

$$\sum_{i=n+1}^{2n+1}\sum_{j=1}^n
(-1)^{i+j+n}(id\wedge\omega)(u_j\cdot \omega(u_1,\ldots,\hat{u_j},\ldots,
u_n,u_i),u_{n+1},\ldots,\hat{u_i},\ldots,u_{2n+1})+$$
$$\sum_{i=n+1}^{2n+1}
(-1)^{i+1}(id\wedge\omega)(u_i\cdot  
\omega(u_1,\ldots,u_n),u_{n+1},\ldots,\hat{u_i},\ldots,u_{2n+1})=$$

$$B_1+B_2+B_3+B_4,$$
where

$$B_1=\sum_{i=n+1}^{2n+1}\sum_{j=1}^n
(-1)^{i+j+n}u_j\cdot \omega(u_1,\ldots,\hat{u_j},\ldots,
u_n,u_i)\cdot \omega(u_{n+1},\ldots,\hat{u_i},\ldots,u_{2n+1}),$$

$$B_2=\sum_{i=n+1}^{2n+1}\sum_{j=1}^n\sum_{n<s\le 2n+1,s\ne i}
(-1)^{i+j+s+1+\delta(s>i)}\cdot$$
$$u_s\cdot \omega(u_j\cdot \omega(u_1,\ldots,\hat{u_j},\ldots,
u_n,u_i),u_{n+1},\ldots,\hat{u_s},\ldots,\hat{u_i},\ldots,u_{2n+1}),$$

$$B_3=\omega(u_1,\ldots,u_n)\cdot
\sum_{i=n+1}^{2n+1}(-1)^{i+1}u_i\cdot  
\omega(u_{n+1},\ldots,\hat{u_i},\ldots,u_{2n+1}),$$

$$B_4=\sum_{i=n+1}^{2n+1}\sum_{n<s\le 2n+1,s\ne i} 
(-1)^{i+s+n+\delta(s>i)}\cdot$$
$$
u_s\cdot\omega(u_i\cdot\omega(u_1,\ldots,u_n),u_{n+1},
\ldots,\hat{u_s},\ldots,
\hat{u_i},\ldots,u_{2n+1}).$$
Notice that 
$$A_1=B_1,$$
$$A_3=B_3.$$

So, to check $(n+1)$-Lie identity for $(n+1)$-multiplication $\tilde\omega$  
we must prove, that 
$$A_2=B_2+B_4.$$
We have 
$$B_2=B_{2,1}+B_{2,2},$$
$$B_4=B_{4,1}+B_{4,2},$$
where 
$$B_{2,1}=
\sum_{i=n+1}^{2n+1}\sum_{j=1}^n\sum_{n<s\le 2n+1,s\ne i}
(-1)^{i+j+s+1+\delta(s>i)}\cdot$$
$$u_s\cdot u_j\cdot \omega( \omega(u_1,\ldots,\hat{u_j},\ldots,
u_n,u_i),u_{n+1},\ldots,\hat{u_s},\ldots,\hat{u_i},\ldots,u_{2n+1}),$$

$$B_{2,2}=
\sum_{i=n+1}^{2n+1}\sum_{j=1}^n\sum_{n<s\le 2n+1,s\ne i}
(-1)^{i+j+s+1+\delta(s>i)}\cdot$$
$$u_s\cdot \omega(u_1,\ldots,\hat{u_j},\ldots,u_n,u_i)\cdot 
\omega(u_j,u_{n+1},\ldots,\hat{u_s},\ldots,\hat{u_i},\ldots,u_{2n+1}),$$

$$B_{4,1}=\sum_{i=n+1}^{2n+1}\sum_{n<s\le 2n+1,s\ne i} 
(-1)^{i+s+n+\delta(s>i)}\cdot$$
$$u_s\cdot u_i\cdot \omega(\omega(u_1,\ldots,u_n),u_{n+1},
\ldots,\hat{u_s},\ldots,
\hat{u_i},\ldots,u_{2n+1}),$$

$$B_{4,2}=\sum_{i=n+1}^{2n+1}\sum_{n<s\le 2n+1,s\ne i} 
(-1)^{i+s+n+\delta(s>i)}\cdot$$
$$u_s\cdot \omega(u_1,\ldots,u_n)\cdot 
\omega(u_i,u_{n+1},\ldots,\hat{u_s},\ldots,
\hat{u_i},\ldots,u_{2n+1}).$$
We see that 
$$B_{4,2}=$$
$$\omega(u_1,\ldots,u_n)\cdot $$
$$
\sum_{i=n+1}^{2n+1}\sum_{n<s\le 2n+1}
(-1)^su_s\cdot \omega(u_{n+1},\ldots,\hat{u_s},\ldots,
u_{2n+1})=$$
$$n\,\omega(u_1,\ldots,u_n)\cdot 
\sum_{n<s\le 2n+1}
(-1)^su_s\cdot \omega(u_{n+1},\ldots,\hat{u_s},\ldots,
u_{2n+1}).$$
Further,
$$B_{2,2}=
\sum_{i=n+1}^{2n+1}\sum_{j=1}^n\sum_{n<s\le 2n+1,s\ne i}(-1)^{i+j+n+1}\cdot $$
$$u_s\cdot \omega(u_1,\ldots,\hat{u_j},\ldots,u_n,u_i)
\cdot 
\omega(u_{n+1},\ldots,
u_{s-1},u_j,u_{s+1},\ldots,\hat{u_i},\ldots,u_{2n+1}).$$
Therefore,
$$B_{2,2}+B_{4,2}=$$
$$\pm n\,\sum_{s=n+1}^{2n+1}u_s\cdot f^{\omega}(u_1,\ldots,u_{n-1},u_n,\ldots,\hat{u_s},\ldots,u_{2n+1})
$$
So, if $n\equiv 0(mod\,p),$ or $f^{\omega}=0,$ then 
$$B_{2,2}+B_{4,2}=0.$$

Let 
$$A_2-B_{2,1}-B_{4,1}=\sum_{1\le i<s\le 2n+1}
u_i\cdot u_s\cdot g_{i,s},$$
where the polynomial $g_{i,s}$ does not depend from $u_i$ and $u_s.$ 
One can check that  $g_{i,s}=0,$ if $1\le i,s\le n$ or 
$n<i,s\le 2n+1$ By $n$-Lie identity for $\omega$ we have 
$g_{i,s}=0,$ if $1\le i\le n, n<s\le 2n+1.$

So, $\tilde\omega$ satisfies $(n+1)$-Lie identity. 

\begin{crl} \label{20augu}  
Let $U$ be an associative commutative algebra with commuting derivations $\der_1,\ldots,\der_n.$ Let $\tilde\omega=id\wedge \der_1\wedge \cdots\wedge \der_n$
Then $(U,\cdot,\tilde\omega)$ is $(n+1)$-Lie with identity $f^{\tilde \omega}=0.$ If $U$ has unit 1, then the algebra $(U,\cdot,\tilde\omega)$ satisfies the identity \rm(\ref{19aug}).
\end{crl}  

The algebra $(U,\cdot,\tilde\omega)$ constructed in corollary 
(\ref{20augu}) is called {\it $n$-Lie Jacobian algebra of type $W.$}

{\bf Proof.} Follows from corollary \ref{20aug} and theorem \ref{14august2}.

\section{\label{poly} Polynomial principle and $\D$-invariants} 

In this section we discuss two methods used in our calculations. 

Call the first method a {\it Polynomial principle}. It means the following. 
Let $U$ be some associative commutative algebra over a field $K$ 
with binary multiplication $U\times U\rightarrow U, (u,v)\mapsto u\cdot v.$ 
Let $\der_1,\ldots,\der_n$ be commuting derivations of $U:$ 
$[\der_i,\der_j]=0,$ for any $i,j=1,\ldots,n.$ Let us given 
a polylinear map $\omega: U\times \cdots \times U\rightarrow U$ with 
$n$ arguments, such that $\omega(u_1,\ldots,u_n)$ is a linear combination of 
elements of the form $\der_1^{i_1}(u_1)\cdot \cdots \cdot \der_n^{i_n}u_n$ 
for any $u_1,\ldots,u_n\in U.$ 

Suppose that some statement concerning associative commutative algebra 
$U,$ $n$-ary polylinear map $\omega: U\times \cdots\times U\rightarrow U$ 
and derivations $\der_1,\ldots,\der_n$  was obtained by using: 
\begin{itemize}
\item Leibniz rule:
$$\der_i(u_1\cdot u_2)=\der_i(u_1)\cdot u_2+u_1\cdot \der_i(u_2),$$
for any $u_1,u_2\in U$ and $i=1,\ldots,n.$   
\item linear properties of $U$  
\item associativity and commutativity properties of $U.$ 
\end{itemize}
Then this statement is true for any associative commutative algebra $\tilde U$ with 
commuting derivations $\tilde\der_1,\ldots,\tilde\der_n.$
In particular, this statement is true for the polynomial algebra 
$U=K[x_1,\ldots,x_n]$ with derivations $\der_i=\der/\der\,x_i, i=1,\ldots,n.$ 

The second method is based on the study of $\mathcal D$-invariants 
\cite{Dzhumavestnik}. Let ${\mathcal D}$ be an abelian subalgebra of 
$Der\,U$ generated by commuting derivations $\der_1,\ldots,\der_n.$ 
Let $C^k(U,U)=\{\psi: U\times\cdots\times U\rightarrow U\}$ 
be the space of skew-symmetric polylinear maps with $k$ arguments 
and $C^*(U,U)=\oplus_kC^k(U,U).$
  
Let $\rho : {\D}\rightarrow End\, C^*(U,U)$ be a representation of an 
abelian Lie algebra $\D$ defined by 
$$\rho(X)\psi(u_1,\ldots,u_k)=$$
$$X (\psi(u_1,\ldots,u_k))-
\sum_{l=1}^n\psi(u_1,\ldots,u_{l-1},X(u_l),u_{l+1},\ldots,u_k),$$
for $\psi\in C^k(U,U).$
We say that $\psi\in C^*(U,U)$ is {\it $\D$-invariant,} if 
$\rho(X)\psi=0,$ for any $X\in \D.$ 

Define $\wedge$ and $\wedge'$ products on $C^*(U,U).$ 
For 
$\psi\in C^k(U,U)$ and 
$\phi\in C^l(U,U)$ set 
$$\psi\smile \phi(u_1,\ldots,u_{k+l})=
\psi(u_{1},\ldots,u_{k})\cdot \phi(u_{k+1},\ldots,u_{k+l}),$$
$$\psi\smile' \phi(u_1,\ldots,u_{k+l})=
\psi(\phi(u_{1},\ldots,u_{l}),u_{l+1},\ldots,u_{k+l}).$$
Define $\psi\wedge \phi\in C^{k+l}(U,U)$ and 
$\psi\wedge' \phi\in C^{k+l-1}(U,U)$ by
$$\psi\wedge \phi(u_1,\ldots,u_{k+l})=$$
$$\sum_{\sigma\in Sym_{k,l}}
sign\,\sigma\,(\psi\smile \phi)(u_{\sigma(1)},\ldots,u_{\sigma(k)},
u_{\sigma(k+1)},\ldots,u_{\sigma(k+l)}),$$
$$\psi\wedge' \phi(u_1,\ldots,u_{k+l})=$$
$$\sum_{\sigma\in Sym_{l,k-1}}
sign\,\sigma\,(\psi\smile'\phi)(u_{\sigma(1)},\ldots,u_{\sigma(l)}, 
u_{\sigma(l+1)},\ldots,u_{\sigma(k+l)}).$$

\begin{prp} \label{9aug} If $\psi, \phi\in C^*(U,U),$ then 
$$\rho(X)(\psi+\phi)=\rho(X)\psi+\rho(X)\phi,$$
$$\rho(X)(\lambda\psi)=\lambda\rho(X)\psi,$$
$$\rho(X)(\psi\wedge \phi)=\rho(X)\psi\wedge \phi+\psi\wedge \rho(X)\phi,$$
$$\rho(X)(\psi \wedge' \phi)=\rho(X)\psi\wedge' \phi+\psi\wedge' \rho(X)\phi.$$
\end{prp}

\begin{crl}\label{22sept} The subspace of $\D$-invariants  
$C^k(U,U)^{\D}=\{\psi\in C^k(U,U): \rho(X)\psi=0, \forall X\in \D\}$ 
is close under $\wedge$ and $\wedge'$ products. 
\end{crl}

So, any polylinear map constructed by $\D$-invariant maps 
using $\wedge$ and $\wedge'$ products will be also $\D$-invariant. 
For example, $Jac_n^S\in C^n(U,U)^{\D}$ and $Jac_{n+1}^W\in 
C^{n+1}(U,U)^{\D}.$ 

Let $pr: U\rightarrow K$ be a projection map: 
$$pr\,\sum_{\alpha\in {\bf Z}^n}\lambda_{\alpha}x^{\alpha}=\lambda_0.$$
Prolong this map to 
$$pr : C^*(U,U)\rightarrow C^*(U,K),$$
$$(pr\,\psi)(u_1,\ldots,u_n)= pr(\psi(u_1,\ldots,u_n)).$$

\begin{thm}\label{altynbala} Let $\psi\in C^k(U,U)^{\D}.$ 
Then 
$$\psi(u_1,\ldots,u_k)=\sum_{\alpha(1),\ldots,\alpha(n)\in{\bf Z}^n_+}
\frac{\der^{\alpha(1)}(u_1)}{\alpha(1)!}\cdots 
\frac{\der^{\alpha(k)}(u_k)}{\alpha(k)!} pr\,\psi(x^{\alpha(1)},\ldots,
x^{\alpha(k)}).$$
\end{thm}

This theorem follows from the results of \cite{Dzhumavestnik}. Call 
$k$-typles $(x^{\alpha(1)},\ldots,x^{\alpha(k)}),$ or simply  
$(\alpha(1),\ldots,\alpha(k)),$ such that 
$pr\,\psi(x^{\alpha(1)},\ldots,x^{\alpha(k)})\ne 0,$ as a {\it support} of 
$\psi.$ So, by theorem \ref{altynbala} to prove $\psi=0$ it is enough to 
establish that $\psi$ does not have any support $k$-typle.

\section{\label{SS} Derivations of Jacobian algebras of type $S$}
 
In this section $char \,K=0$ and 
$U=K_n$ or $K_n^+$ and $n>1,$ if otherwise is not stated. 
Elements of $U$ are denoted as 
$u,v,w, u_1,v_1,w_1,\ldots .$
Let $\theta=\sum_{i=1}^n\epsilon_i,$ where $\epsilon_i=(0,\ldots,0,\mathop{1}\limits_i,0,\ldots,0)\in {\bf Z}^n.$ 

There are four Cartan Type Lie algebras of formal vector fields. 
One of them is called Special. This algebra is defined as an algebra of 
divergenceless vector fields
$$S_{n-1}=S_{n-1}(U)=<X=u_i\der_i: Div\,X=\sum_{i=1}^n\der_i(u_i)=0>.$$
Recall that $S_n(K_n^+)$ is simple, but $L=S_{n-1}(K_n)$ is not simple: 
its commutant $[L,L]$ is generated by derivations 
$D_{ij}=\der_i(u)\der_j-\der_j(u)\der_i, u\in U$ and 
$$L/[L,L]=<D_{ij}x^{-\theta}, x^{-\theta+\epsilon_i}\der_i: i=1,\ldots,n>\cong 
K^{n+1}.$$

\begin{thm}\label{S} Let $U=K_n.$
 
{\rm i)} The algebra 
$(U,Jac_{n}^S)$ is not simple. It has $1$-dimensional center $<1>$ and the ideal of codimension 1: 
$\bar U =< x^{\alpha}: \alpha\ne -\theta>.$ The factor-algebra 
$(\bar U/<1>, Jac_n^S)$ is simple $n$-Lie algebra. 

{\rm ii)} Takes place the exact sequence 
$$0\rightarrow Out\,(U,\cdot,Jac_{n}^S)\rightarrow Out\,(U,Jac_n^S)
\rightarrow K^2\rightarrow 0.$$
More exactly, the factor-algebra 
$ Out\,(U,Jac_n^S)/Out\,(U,\cdot,Jac_{n}^S)$
is $2$-dimensional and is generated by derivations $\Delta$ 
and $D_{-\theta},$ defined by 
$$\Delta=\sum_{i=1}^nx_i\der_i+n(1-n)^{-1}: x^{\alpha}\mapsto (|\alpha|+n(1-n)^{-1})x^{\alpha},$$
$$D_{-\theta}: x^{\alpha}\mapsto \delta_{\alpha,-\theta}.$$

{\rm iii)} The linear maps $D_i:U\rightarrow U, 
i=1,\ldots,n,$  given by 
$$D_i: u\mapsto x^{-\theta+\epsilon_i}\der_i(u),$$
are outer derivations of $(U,\cdot,Jac_n^S).$ Their classes of derivations 
in $Out\,(U,\cdot,Jac_n^S)=
Der\,(U,Jac_n^S)/Int\,(U,Jac_n^S)$ form a base. 
In particular, $dim\,Out\,(U,\cdot,Jac_n^S)=n$ and 
$dim\,Out\,(U,Jac_n^S)=n+2.$

{\rm iv)} Classes of derivations 
$\Delta, D_1,\ldots,D_n$ in $Out\,(\bar U, Jac_n^S)/<1>)$ form a base. 
In particular, $dim\,Out\,(\bar U, Jac_n^S)/<1>=n+1.$ 
\end{thm}

\begin{crl}
Let $U=K[x_1^{\pm \frac{1}{n-1}},\ldots,x_n^{\pm \frac{1}{n-1}}].$ 
Then eigenspaces of the derivation $\Delta$ endow the algebra 
$(U,Jac_n^S)$ by grading. If 
$$U_{[k]}=
\{u\in U:\Delta(u)=k\,u\}=\{x^{\alpha} : |\alpha|=k+n(n-1)^{-1}\},$$  
then 
$$u_1\in U_{[k_1]},\ldots, u_n\in U_{[k_n]}\Rightarrow  Jac_n^S(u_1,\ldots,u_n)
\in U_{[k_1+\cdots+k_n]}.$$ 
Here  $k\in (n-1)^{-1}{\bf Z}.$ 
\end{crl}

\begin{thm}\label{S^+} Let $U=K_n^+.$ Then 
the algebra $(U^+,Jac_n^S)$ has 
$1$-dimensional center $<1>.$ Its factor-algebra $(U^+/<1>,Jac_n^S)$ is 
simple $n\mbox{-Lie}.$ It has one outer derivation $\Delta.$ In particular, 
$dim\,Out\, (U,Jac_n^S)=$ \\ $dim\,Out\, (U/<1>,Jac_n^S)=1.$ 
\end{thm}

Theorem~\ref{S^+} follows from theorem~\ref{S}. Therefore, below we 
can consider only the case $U=K_n.$ 

\begin{lm} \label{3au} 
Any divergenceless derivation of $U$ is a linear combination of derivations 
$D_{i,j}(u)$ and derivations $x^{-\theta+\epsilon_i}\der_i.$ 
\end{lm}

\begin{lm}\label{3aug} Any derivation of the form $D_{i,j}(u),$ where $i<j, i,j =1,\ldots,n, u\in U,$ can be presented as an interior derivation 
$L_{u_1,\ldots,u_{n-1}}$ of the Jacobian algebra $Jac_n^S.$ 
\end{lm}

{\bf Proof.}  It is easy to see that
$$Jac_n^S(x_1,\ldots,x_{i-1},u,x_{i+1},\ldots,x_{j-1},v,x_{j+1},\ldots,x_n)=
$$
$$\der_i(u)\der_j(v)-\der_j(u)\der_i(v).$$
Therefore,
$$D_{i,j}(u)=(-1)^{n-j}ad\,(x_1,\ldots,x_{i-1},u,x_{i+1},\ldots,x_{j-1},x_{j+1},\ldots,x_n).$$

\begin{lm}\label{3augu} The element $x^{-\theta}$ cannot be presented as a 
Jacobian $Jac_n^S(u_1,\ldots,u_n).$ 
Any element $x^{\beta},$ where $\beta\in {\bf Z}^n, 
\beta\ne -\theta,$ can be presented as a 
linear combination of elements $Jac_n^S(u_1,\ldots,u_n),$ for some 
$u_1,\ldots,u_n\in U.$ 
\end{lm}

In other words, the first cohomology group of $(U,Jac_n^S)$ with coefficients 
in the trivial module is 1-dimensional. 

{\bf Proof.} Since $Jac_n^S: \wedge^n U\rightarrow U$ is polylinear, 
if $x^{-\theta}=Jac_n^S(u_1,\ldots,u_n),$ for some $u_1,\ldots,u_n\in U,$ then 
$$\lambda x^{-\theta}= Jac_n^S(x^{\alpha_1},\ldots,x^{\alpha_n}),$$
for some $\lambda\in K, \alpha_1,\ldots,\alpha_n\in {\bf Z}^n.$ 
Thus,
$$\lambda=x^{\theta} Jac_n^S(x^{\alpha_1},\ldots,x^{\alpha_n})=$$
$$\sum_{\sigma\in Sym_n}sign\,\sigma\,x_{\sigma(1)}\der_{\sigma(1)}
(x^{\alpha_1})\cdots x_{\sigma(n)}\der_{\sigma(n)}(x^{\alpha_n})=$$
$$(\sum_{\sigma\in Sym_n}sign\,\sigma\,\alpha_{1,\sigma(1)}\cdots \alpha_{n,\sigma(n)})x^{\alpha_1+\cdots+\alpha_n},$$
where $\alpha_i=(\alpha_{i,1},\ldots,\alpha_{i,n})\in {\bf Z}^n.$ 
Therefore, we come to two conclusions: the first  
$$\alpha_1+\cdots+\alpha_n=0\in{\bf Z}^n,$$
and the second,
$$\lambda=det (\alpha_{i,j}).$$
The first relation means that the sum of rows of the matrix $(\alpha_{i,j})$ is $0.$ Therefore, its determinant is also $0$ and by the second relation 
$\lambda=0.$

Let $\beta\ne -\theta.$ Then $\beta_i\ne -1$ for some $i.$ Therefore, 
$\beta_i+1\ne 0,$ and 
$$x^{\beta}=Jac_n^S(x_1,\ldots,x_{i-1},(\beta_i+1)^{-1}
x^{\beta+\epsilon_i},x_{i+1},\ldots,x_n).$$

\begin{crl}\label{22septpopevastane} $D_{-\theta}\in Der\,(U,Jac_n^S).$
\end{crl}

\begin{lm} \label{4aug}
The element $x^{-\theta}$ cannot be presented in the form  
$x^{-\theta+\epsilon_i}\der_i(u),$ where $u\in U, i=1,\ldots,n.$ 
\end{lm}

{\bf Proof.} Suppose that $\lambda x^{-\theta}= 
x^{-\theta+\epsilon_i}\der_i(x^{\alpha}),$ for some 
$\alpha\in {\bf Z}^n$ and $\lambda\in K.$ We have 
$$x^{-\theta}=x^{-\theta+\epsilon_i}\der_i(x^{\alpha})=\alpha_ix^{-\theta+\alpha}.$$
Thus   $\alpha=0,$ and in particular, $\alpha_i=0.$ So, $\lambda=0.$

\begin{lm}\label{3augus} The linear maps $D_i=x^{-\theta+\epsilon_i}\der_i: 
u\mapsto x^{-\theta+\epsilon_i}\der_i(u)$ are derivations of $(U,Jac_n^S).$ 
\end{lm}

{\bf Proof.} Take $T_{i,j}=x^{-\theta+\epsilon_i+\epsilon_j}\ln x_j, i<j.$ 
Then $T_{i,j}\not\in U=K_n,$ but 
$T_{i,j}\in \tilde U,$ where $\tilde U=
K[x_1^{\pm 1},\ldots,x_n^{\pm 1},T_{i,j}]$ is an extension of 
$U$ by $T_{i,j}.$
We see that 
$D_i=-D_{i,j}(T_{i,j}).$ 
So, by lemma \ref{3aug}, 
$D_i$ is a derivation of $Jac_n^S(\tilde U).$ 
In particular, $D_i$ satisfies the Leibniz rule 
$$D_i Jac_n^S(u_1,\ldots,u_n)=\sum_{l=1}^n
Jac_n^S(u_1,\ldots,u_{l-1},D_i(u_l),u_{l+1},\ldots,u_n),$$
for elements $u_1,\ldots,u_n$ of the  
subalgebra $(U,Jac_n^S)\subset (\tilde U,Jac_n^S).$ 
Since $D_i(u)\in U,$ for any $u\in U,$ this means that 
$D_i\in Der\,(U,J_n^S).$ 

\begin{lm}\label{3august}
The linear map $\Delta: U\rightarrow U$ given by 
$$\Delta=\sum_{i=1}^nx_i\der_i+n(1-n)^{-1}$$
is a derivation of the algebra $(U,Jac_n^S).$ 
\end{lm}

{\bf Proof.} Define polynomial $f$ by 
$$f(u_1,\ldots,u_n)=\Delta (Jac_n^S(u_1,\ldots,u_n))
-\sum_{l=1}^nJac_n^S(u_1,\ldots,u_{l-1},\Delta(u_l),u_{l+1},\ldots,u_n).$$

Notice that $[\der_i,\Delta]=\der_i.$ Since $Jac_n^S$ is $\D$-invariant, 
$$\rho(\der_i)f=\rho([\der_i,\Delta]) Jac_n^S=\rho(\der_i)Jac_nS=0.$$
So, $f$ is $\D$-invariant and by theorem \ref{altynbala},
$$f(u_1,\ldots,u_n)=
\sum_{\alpha(1),\ldots,\alpha(n)\in {\bf Z}^n_+}
\frac{\der^{\alpha(1)}(u_1)}{\alpha(1)!}\cdots \frac{\der^{\alpha(n)}(u_n)}{\alpha(n)!} pr\,f(x^{\alpha(1)},\ldots,x^{\alpha(n)}),$$
where $|\alpha(1)|+\cdots+|\alpha(n)|=n.$ 

We see that $pr\,f: \wedge^nU\rightarrow K$ is uniquely defined by 
$f(x_1,\ldots,x_n).$ In other words,  
$$Support(f)\subseteq  C^n(R(\pi_1),K)\cong <x_1\wedge \cdots \wedge x_n>\cong K.$$
We have,
$$\Delta(1)=n(1-n)^{-1},\quad  \Delta(x_i)=(1-n)^{-1}, \quad Jac_n^S(x_1,\ldots,x_n)=1.$$
Therefore, 
$$f(x_1,\ldots,x_n)=n(1-n)^{-1}-\sum_{l=1}^n(1-n)^{-1}=0.$$
So, $Support(f)=\emptyset,$ and $pr\,f=0.$ 
Thus, by theorem \ref{altynbala}, $f=0.$  
In other words, $\Delta\in Der\,(U, Jac_n^S).$

\begin{lm}\label{outer} Derivations $\Delta, D_1,\ldots,D_n\in Der\,
(U,Jac_n^S)$ 
are outer. Moreover, any non-trivial linear combination of derivations 
$D_1,\ldots,D_n$ and $\Delta$ is outer.
\end{lm}

{\bf Proof.}  Suppose that 
$$\lambda_0\Delta(v)+\sum_{i=1}^n\lambda_iD_i(v)=
Jac_n^S(u_1,\ldots,u_{n-1},v),$$
for some $\lambda_0,\lambda_1,\ldots,\lambda_n\in K, u_1,\ldots,u_{n-1}\in U$ 
and for any  $v\in U.$ 

For $v=1$ we obtain that $n(1-n)^{-1}\lambda_0=0,$ and $\lambda_0=0.$ 

For $v=x_i$ we have 
$$\lambda_ix^{-\theta+\epsilon_i}=(-1)^{n-i}Jac_{n-1}^S(u_1,\ldots,u_{n-1}).$$
By lemma \ref{3augu}, this is possible only for the case of $\lambda_i=0.$  

\begin{lm} \label{6aug} Suppose that $D_{i,j}(x^{\beta})u=0, u\in U,$ 
for any $\beta\in {\bf Z}^n_+, |\beta|>2.$ Then $u\in <1>.$ 
\end{lm}

{\bf Proof.} Suppose that $u=\sum_{\alpha\in {\bf Z}^n, q\le |\alpha|\le s}
\lambda_{\alpha}x^{\alpha}$ and $\lambda_{\bar\alpha}\ne 0,$ for some 
$\bar\alpha\in {\bf Z}^n,$ such that $|\bar\alpha|=s$ and $\bar\alpha\ne 0.$ 
Then $\bar\alpha_i\ne 0,$ for some $i.$ 
Take $\beta=3\epsilon_j,$ for some $j\ne i.$ 
Recall that we suppose that $n>1.$
Then $D_{i,j}(x^{\beta})= -3x_j^2\der_i$ and 
$$D_{i,j}(x^{\beta})u=
-3\lambda_{\bar\alpha}\bar\alpha_i x^{\bar\alpha+2\epsilon_j}+ w,$$ 
for some $w\in U.$ Notice that $w$ is a linear combination of $x^{\gamma},$ 
where $\gamma\in {\bf Z}^n, |\gamma|\le s+2$ and $\gamma\ne \bar\alpha+2\epsilon_j.$   
Therefore, $-3\lambda_{\bar\alpha}\bar\alpha_i=0,$ i.e., $\lambda_{\bar\alpha_i}=0.$
The obtained contradiction shows that $s=0.$ 
By similar arguments it is easy to obtain that $q=0.$

\begin{lm}\label{6august} \label{U^+} Let $D\in Der\,
(U, Jac_n^S),$ and $D(u)=0,$ for any 
$u\in U^+.$ Then $D=\lambda D_{-\theta},$ for some $\lambda\in K.$   
\end{lm}

{\bf Proof.} Suppose that $D(x^{\alpha})\ne 0,$ for some $\alpha\in {\bf Z}^n, {\alpha}\not\in {\bf Z}^n_+$ and $D(x^{\beta})=0,$ for any $\beta\in 
{\bf Z}^n, |\beta|>|\alpha|.$ 
Notice that 
$$Jac_n^S(x_1,\ldots,x_{i-1},u,x_{i+1},\ldots,x_{j-1},v,x_{j+1},\ldots,x_n)=
D_{i,j}(u)v=-D_{i,j}(v)u.$$
For any $\beta\in {\bf Z}^n, |\beta|>2,$   
$$
|Jac_n^S(x_1,\ldots,x_{i-1},x^{\alpha},x_{i+1},\ldots,x_{j-1}, x^{\beta}
,x_{j+1},\ldots,x_n)|> |x^{\alpha}|.$$
So,  
$$D(Jac_n^S(x_1,\ldots,x_{i-1},x^{\alpha},x_{i+1},\ldots,x_{j-1}, x^{\beta}
,x_{j+1},\ldots,x_n))=0,$$
if $\beta\in {\bf Z}^n_+, |\beta|>2.$ Since $D$ is derivation, and 
$D(u)=0,$ for any $u\in U^+,$ 
$$D(Jac_n^S(x_1,\ldots,x_{i-1},x^{\alpha},x_{i+1},\ldots,x_{j-1},x^{\beta}
,x_{j+1},\ldots,x_n))=$$
$$Jac_n^S
(x_1,\ldots,x_{i-1},D(x^{\alpha}),x_{i+1},\ldots,x_{j-1}, x^{\beta}
,x_{j+1},\ldots,x_n)),$$
Thus,
$$D_{i,j}(x^{\beta})(D(x^{\alpha}))=0,$$
for any $\beta\in {\bf Z}^n_+, |\beta|>2.$ 
By lemma~\ref{6aug}, 
$$D(x^{\alpha})\in \cap_{s=1}^nKer\der_s=<1>.$$
So, 
\begin{equation}\label{6augu}
D(x^{\alpha})=\lambda \cdot 1,
\end{equation}
for some $\lambda\in K.$  

Suppose that $\lambda\ne 0.$

Now take $\beta=\epsilon_j.$ Then by the Leibniz rule 
$$D(\der_i(x^{\alpha}))=D(Jac_n^S(x_1,\ldots,x_{i-1},x^{\alpha},x_{i+1},\ldots,x_n))=$$
$$Jac_n^S
(x_1,\ldots,x_{i-1},D(x^{\alpha}),x_{i+1},\ldots,x_n))=\der_i(D(x^{\alpha})),$$
and, according to (\ref{6augu}),
\begin{equation}\label{6augus}
D(\der_i(x^{\alpha}))=0,
\end{equation}
for any $i=1,\ldots,n.$ 

We have
$$Jac_n^S(x_1,\ldots,x_{i-1},x^{\alpha},x_{i+1},\ldots,x_{j-1},
x_ix_j,x_{j+1},\ldots,x_n)=$$
$$x_i\der_i(x^{\alpha})-x_j\der_j(x^{\alpha}).$$
Therefore, by (\ref{6augu}) and the  Leibniz rule, 
$$(\alpha_i-\alpha_j)D(x^{\alpha})=$$
$$D(Jac_n^S(x_1,\ldots,x_{i-1},x^{\alpha},x_{i+1},\ldots,x_{j-1},
x_ix_j,x_{j+1},\ldots,x_n))=0.$$
So, if $\lambda\ne 0,$ then 
\begin{equation}
\label{7aug}\alpha_i=t,
\end{equation} for some $0\ne t\in {\bf Z}$ and for any $i=1,\ldots,n.$ 

Take  $\beta=\epsilon_i+2\epsilon_j.$ 
We have 
$$Jac_n^S(x_1,\ldots,x_{i-1},\der_j(x^{\alpha}),
x_{i+1},\ldots,x_{j-1},x_ix_j^2, x_{j+1},\ldots,x_n)=
(2\alpha_i-\alpha_j+1)x^{\alpha}.$$
Then by (\ref{6augu}), (\ref{7aug}) we obtain that  
$$(2\alpha_i-\alpha_j+1) D(x^{\alpha})=$$
$$D(Jac_n^S(x_1,\ldots,x_{i-1},\der_j(x^{\alpha}),
x_{i+1},\ldots,x_{j-1},x_ix_j^2, x_{j+1},\ldots,x_n))=0.$$
So, by (\ref{7aug}), $t=-1$
and $D=\lambda D_{-\theta}.$

\begin{lm} \label{combin} 
Any derivation of $(U,Jac_n^S)$ is a linear combination 
of some interior derivation and derivations $\Delta, D_{-\theta}, D_1,\ldots,D_n.$
\end{lm}

{\bf Proof.} By lemma \ref{6august} any derivation $D$ is defined by 
its restriction on $U^+$ up to $D_{-\theta}.$  
Any linear map $F: U^+\rightarrow U$ is 
a linear combination of linear maps of the form $x^{\alpha}\der^{\beta}: 
U\rightarrow U, v\mapsto x^{\alpha}(\der^{\beta}(v)),$ where $\alpha 
\in {\bf Z}^n, \beta\in {\bf Z}^n_+.$ 
Suppose that $F=\sum_{\alpha\in {\bf Z}^n,\beta\in {\bf Z}^n_+}
\lambda_{\alpha,\beta}x^{\alpha}\der^{\beta}\in Der\,(U,J_n^S).$
%We will prove that $F$ can be presented as a linar combination of derivations 

{\bf Step 1}. Any derivation $F$ is a linear combination of $\Delta$ and 
some derivation $F'$ with the property $F'(1)=0.$  Let us prove it.   
Since $1$ is central, $F(1)$ is also central for any derivation $F.$
The center of $U$ is 1-dimensional and is generated by one element $1.$ 
Therefore, $F(1)=\lambda_{0,0}.$
If $\lambda_{0,0}\ne 0,$ we can take instead of $F$ 
$$F'=F-{\lambda_{0,0}(1-n)}{n}^{-1}\Delta.$$
Then $F'\in Der\,(U,Jac_n^S)$ and $\lambda'_{0,0}=0,$ if 
$$F'=\sum_{\alpha\in{\bf Z}^n, \beta\in {\bf Z}^n_+}
\lambda'_{\alpha,\beta}x^{\alpha}\der^{\beta}.$$

{\bf Step 2}. Suppose that the derivation $F'$ satisfies the following 
normalisation
\begin{equation} \label{norm1} F'(1)=0 \end{equation}
Prove that $F'$ can be presented as a linear combination of 
derivations $D_1,\ldots,D_n$ and some interior derivation and some 
differential operator of differential order $>1.$   

Take 
$$\bar F'=F'-\sum_{i=1}^n\lambda'_{-\theta+\epsilon_i,\epsilon_i} D_i.$$
Then $F''\in Der\,(U,Jac_n^S)$ and 
$$\bar F'(x_i)\in <x^{\alpha}: \alpha\ne -\theta+\epsilon_i>,$$
for any $i=1,\ldots,n.$ 
Let $\bar F'=\sum_{\alpha\in {\bf Z}^n,\beta\in {\bf Z}^n_+}
\bar \lambda'_{\alpha,\beta} x^{\alpha}\der^{\beta}.$ Then 
$$\bar\lambda_{\alpha,\epsilon_1}\ne 0\Rightarrow \alpha_j\ne -1, \mbox{\;for some \;} j.$$
Take 
$$F''=\bar F'+\bar\lambda_{\alpha,\epsilon_1}(\alpha_j+1)^{-1}
D_{1,j}(x^{\alpha+\epsilon_j}).$$
It is easy to see that 
$F''\in Der\,(U,Jac_n^S),$ and 
$$\lambda''_{\alpha,\epsilon_1}=0,$$
where
$F''=\sum_{\alpha\in {\bf Z}^n,\beta\in {\bf Z}^n_+} 
\lambda''_{\alpha,\beta}x^{\alpha}\der^{\beta}.$ 

Let us use induction on $l$ and suppose that $F''$ some $l>1$ has the 
following property:  
$\lambda''_{\alpha,\epsilon_i}=0,$ for any $\alpha\in {\bf Z}^n,$  
if $i<l$ and  $\lambda''_{\alpha(l),\epsilon_l}\ne 0$ for some $\alpha(l)\in 
{\bf Z}^n.$ 

Let $(\tilde U,Jac_{n-l+1}^S)$ be Jacobian algebra of type $S$ with 
vector space $\tilde U=K[x_l^{\pm 1},\ldots,x_n^{\pm 1}]$ and 
multiplication $Jac_{n-l+1}^S=\der_l\wedge\cdots
\wedge \der_n.$ 

Let $G$ be the restriction of $F''$ to subalgebra 
$(\tilde U,Jac_{n-l+1}^S).$ Then 
$$G(u)=\sum_{\alpha'=(\alpha_1,\ldots,\alpha_{l-1})\in {\bf Z}^{l-1}}
x^{\alpha'} G_{\alpha'}(u),$$ 
for any $u\in \tilde U,$ where $G_{\alpha'}$ is a linear map 
defined on $\tilde U$ with coefficients in $\tilde U.$ 
  
Since $F''(x_i)=0, i<l$ and $Jac_{n-l+1}^S=i(x_1)\cdots i(x_{l-1})Jac_n^S,$  
then 
$$G(Jac_{n-l+1}^S(u_l,\ldots,u_n))=$$
$$F''Jac_n^S(x_1,\ldots,x_{l-1},u_l,\ldots,u_n)=$$
$$\sum_{s=l}^n
Jac_n^S(x_1,\ldots,x_{l-1}, u_l,\ldots,u_{s-1},F''(u_s),u_{s+1},\ldots,u_n)
=$$ 
$$\sum_{s=l}^nJac_{n-l+1}(u_l,\ldots,u_{s-1},G(u_s),u_{s+1},\ldots,u_n),$$
for any $u_l,\ldots,u_n\in <x^{\alpha}: \alpha\in{\bf Z}^n, 
\alpha_i=0, i<l >.$ 
Therefore,
$$\sum_{\alpha'\in {\bf Z}^{l-1}}x^{\alpha'}G_{\alpha'}
(Jac_{n-l+1}(u_l,\ldots,u_n))=$$
$$\sum_{s=1}^n\sum_{\alpha'\in {\bf Z}^{l-1}}x^{\alpha'} 
Jac_{n-l+1}(u_l,\ldots,u_{s-1},G_{\alpha'}(u_s),u_{s+1},\ldots,u_n)=$$
$$\sum_{\alpha'\in {\bf Z}^{l-1}}x^{\alpha'} \sum_{s=1}^n
Jac_{n-l+1}(u_l,\ldots,u_{s-1},G_{\alpha'}(u_s),u_{s+1},\ldots,u_n).$$
In other words, $G_{\alpha'}\in Der\,(\tilde U,J_{n-l+1}^S),$ 
for any $\alpha'\in 
{\bf Z}^{n}, \alpha'_l=\ldots=\alpha'_n=0.$ 
By an inductive suggestion $G_{\alpha'}$ is a linear combination 
of derivations of $(\tilde U,Jac_{n-l+1}^S$ of the form 
$x^{-\sum_{s=l}^n\epsilon_s+\epsilon_j}\der_j,$ 
where $j=l,\ldots,n$ and a derivation of the form $L_{u_l,\ldots,u_{n-1}},$ 
for some $u_l,\ldots,u_{n-1}\in \tilde U.$ Here the derivation $\tilde \Delta=
\sum_{s=1}^nx_s\der_s+(n-l+1)/(l-n)$ does not appear because of normalisation 
(\ref{norm1}).  Thus,
$x^{\alpha'}G_{\alpha'}$ is a linear combination of derivations 
$A_i:=x^{\alpha'-\sum_{s=1}^n\epsilon_s+\epsilon_i}\der_i, i\ge l,$ and 
some linear operator $B_{\alpha'}:=x^{\alpha'}L_{u_l,\ldots,u_{n-1}}.$ 
Notice that $Div\,A_i=0, i\ge l.$ Moreover, $A_i\in Int\,(U,Jac_n^S),$ because 
of $\lambda''_{-\theta+\epsilon_i,\epsilon_i}=0.$   

By an inductive suggestion and by lemma \ref{3au}, the divergenceless derivation  $L_{u_l,\ldots,u_{n-1}}$ is a linear combination of derivations of the form 
$D_{i,j}(u), l\le i<j, u\in\tilde U.$ Therefore, $B_{\alpha'},$ where 
$\alpha'=(\alpha_1,\ldots,\alpha_{l-1},0,\ldots,0)\in {\bf Z}^n,$ is 
also a linear combination of derivations $D_{i,j}(x^{\alpha'}u), l\le i<j.$  
In particular, $B_{\alpha'}\in Int\,(U,Jac_n^S).$ 

So, we have proved that  $G$ is a linear combination of 
interior derivations. Therefore, we can add to $F''$ some linear 
combination of interior derivations and obtain a new derivation 
$F'''=\sum_{\alpha\in {\bf Z}^n,\beta\in {\bf Z}^n_+} 
\lambda'''_{\alpha,\beta}x^{\alpha}\der^{\beta},$ 
such that 
$\lambda'''_{\alpha,\epsilon_i}=0,$ for any $\alpha\in {\bf Z}^n,$  
if $i<l+1.$ 

Thus, induction on $l$ is possible and finally we obtain that 
$F''$ is a sum of some interior derivation and some linear operator 
$$H= \sum_{\alpha\in {\bf Z}^n,\beta\in {\bf Z}^n_+, |\beta|>1} 
\mu_{\alpha,\beta}x^{\alpha}\der^{\beta}.$$

{\bf Step 3.} Prove that $H=0.$ Present $H$ in the form 
$H=\sum_{\beta\in {\bf Z}^n_+, |\beta|>1} H_{\beta}\der^{\beta},$
where $H_{\beta}=\sum_{\alpha\in {\bf Z}^n}\lambda_{\alpha,\beta}x^{\alpha}.$ 
Suppose that $H\ne 0$ and $H_{\bar\beta}\ne 0,$ for some $\bar\beta\in {\bf Z}^n_+, |\bar\beta|=t>1$ and $H_{\beta}=0,$ for any $\beta\in {\bf Z}^n_+, |\beta|<t.$ We have 
$$Jac_n^S(x_1,\ldots,x_{i-1},x^{\bar\beta},x_{i+1},\ldots,x_n)=
\der_i(x^{\bar\beta}),$$
$$H(x_i)=0, \forall i=1,\ldots,n,$$
Thus,
$$H(Jac_n^S(x_1,\ldots,x_{i-1},x^{\bar\beta},x_{i+1},\ldots,x_n))=
H(\der_i(x^{\bar\beta})=0.$$
On the other hand, by the Leibniz rule 
$$H(Jac_n^S(x_1,\ldots,x_{i-1},x^{\bar\beta},x_{i+1},\ldots,x_n))=$$
$$Jac_n^S(x_1,\ldots,x_{i-1},H(x^{\bar\beta}),x_{i+1},\ldots,x_n)=
\der_i(H(x^{\bar\beta})).$$
So, $H(x^{\bar\beta})\in \cap_{i=1}^nKer\,\der_i=<1>.$ 
Here $\bar\beta$ is any element of ${\bf Z}^n_+,$ such that 
$H_{\bar\beta}\ne 0, |\bar\beta|>1.$ In particular, 
$H_{\epsilon_i+\epsilon_j}=0$ or $H_{\epsilon_i+\epsilon_j}\in <1>.$ 
Therefore, in any case by the Leibniz rule 
$$(\bar\beta_i-\bar\beta_j)H(x^{\bar\beta})=$$
$$H(Jac_n^S(x_1,\ldots,x_{i-1},x^{\bar\beta},x_{i+1},
\ldots,x_{j-1},x_ix_j,x_{j+1},\ldots,x_n))=0.$$
So, $\bar\beta_i=q, i=1,\ldots,n,$ for some $0<q\in {\bf Z}.$ 
In other words, $\bar\beta=q\theta.$    
In particular, 
$t=|\bar\beta|=nq\ge n,$
and, since $n>1,$
$t=|\bar\beta|=nq> q.$
Thus
$q+1\le t,$ and $H(x_i^{q+1})=0$ or $H(x_i^{q+1}\in <1>.$ We have 
$J_n^S(x_1^{q+1},\ldots,x_n^{q+1})=(q+1)^nx^{q\theta},$ 
So, by the Leibniz rule, 
$$(q+1)^nH(x^{q\theta})=$$
$$H(Jac_n^S(x_1^{q+1},\ldots,x_n^{q+1}))=0.$$
Hence, $H_{\bar\beta}=0.$ Contradiction.

{\bf Proof of theorem \ref{S}.} Simplicity of $(U,Jac_n^S)$ was established in 
\cite{Filipov1}. Other statements follow from  
corollary \ref{22septpopevastane} and lemmas \ref{3augus}, \ref{6august} \ref{3august}, \ref{outer} and \ref{combin}.

\section{\label{WW} Derivations of Jacobian algebras of type $W$} 
In this section $K$ is the field of characteristic $0$ 
and  $U=K_n$ or $K_n^+.$
Let 
$$W_n=W_n(U)=<u\der_i: u\in U, i=1,\ldots,n>$$ 
be General Cartan 
Type Lie algebra. This algebra is also called Witt algebra. 

\begin{thm}\label{W}  The $n$-Lie algebra $(U,Jac_{n+1}^W)$ is 
simple. Any derivation of $(U,Jac_{n+1}^W)$ is interior. Moreover, 
$Int\,(U,Jac_{n+1}^W)$ is isomorphic to Witt algebra $W_n.$ 
\end{thm}
 
On $U$ we have two multiplications $Jac_{n+1}^W$ and $Jac_n^S.$ 
They are related by 
$$Jac_{n+1}^W(u_1,\ldots,u_{n+1})=$$
$$(-1)^{n+1}Jac_n^S(u_1,\ldots,u_n)u_{n+1}+
\sum_{i=1}^{n}
(-1)^{i+n-1}Jac_{n,\hat{i}}^W(u_1,\ldots,u_{n})\der_i(u_{n+1}),$$
where
$$Jac_{n,\hat{i}}^W=id\wedge \der_1\wedge\cdots 
\wedge\hat{\der_i}\wedge\cdots\wedge\der_n.$$

Let $W_n+U$ be a semidirect sum:
$$[X+u,Y+v]=[X,Y]+X(v)-Y(u),$$
for any $X,Y\in W_n, u,v\in U.$ 

\begin{lm}\label{07aug} 
The map $W_n\rightarrow W_n+U$ given by 
$X\mapsto X+\lambda Div\,X$ is a monomorphism for any $\lambda\in K.$
\end{lm}

{\bf Proof.} The statement is evident, if $\lambda=0.$
Suppose that $\lambda\ne 0.$ 

The divergence map $X\mapsto Div\,X$ 
is $1$-cocycle in $Z^1(W_n,U):$ 
$$Div[X,Y]=X(Div\,Y)-Y(Div\,X).$$
Therefore, for $f=id+\lambda\,Div,$  
$$f[X,Y]=[X,Y]+\lambda Div[X,Y]=$$
$$[X,Y]+\lambda(X(Div(Y))-Y(Div(X)))=$$
$$[X+\lambda\, Div\,X,Y+\lambda\, Div\,Y]=[f(X),f(Y)].$$
So, $f$ is the homomorphism of Lie algebras. 

Suppose that $X\in Ker\,f.$ Then 
$$X(u)+\lambda (Div\,X) u=0,$$
for any $u\in U.$ Take $u=1.$ Then
$\lambda (Div\,X)=0.$ Thus, $X=\sum_{i=1}^n\mu_i\der_i.$ 
Take now  $u=x_i.$ We have $0=X(u)=\lambda_i.$ So, $X=0.$ 

\begin{lm}\label{9augu} 
For any associative commutative algebra $U$ with commuting derivations 
$\der_1,\ldots,\der_n$ the algebra $(U,Jac_{n+1}^W)$
is $(n+1)$-~Lie and graded:
$$J_{n+1}^W(U_{[i_1]},\ldots,U_{[i_k]})\subseteq U_{[i_1+\ldots+i_k]},$$
where
$$U_{[s]}=<x^{\alpha}: |\alpha|=s+1>, s\in {\bf Z}.$$
\end{lm}

{\bf Proof.} By corollary \ref{20augu} $(U,Jac_{n+1}^W)$ is $(n+1)$-Lie. 
It is graded: 
$$Jac_{n+1}^W(x^{\alpha(0)},\ldots,x^{\alpha(n)})\in 
<x^{\alpha(0)+\cdots+\alpha(n)-\theta}>,$$ 
and
$$|Jac_{n+1}^W(x^{\alpha(0)},\ldots,x^{\alpha(n)})|= 
|\alpha(0)+\cdots+\alpha(n)-\theta|=(|\alpha(0)|-1)+\ldots+(|\alpha(n)|-1)+1.$$ 

\begin{lm}\label{10august} Any interior derivation of $(U,Jac_{n+1}^W)$ 
has the form $X-~{n}^{-1} Div\,X,$ where $X\in W_n(U).$ 
\end{lm}

{\bf Proof.} Any interior derivation  
$$L_{u_1,\ldots,u_n}: v\mapsto Jac_{n+1}^W(u_1,\ldots,u_n,v)$$
can be presented as a sum of some differential operator of first 
order and the operator of multiplication to some element of $U$ 
$$L_{u_1,\ldots,u_n}=X_{u_1,\ldots,u_n}+R_{u_1,\ldots,u_n},$$
where 
$$X_{u_1,\ldots,u_n}=\sum_{i=1}^{n}
(-1)^{n+i-1}Jac_{n,\hat{i}}^W(u_1,\ldots,u_{n})\der_i,$$
$$R_{u_1,\ldots,u_n}v=
(-1)^{n+1}{Jac_{n}^S(u_1,\ldots,u_{n})}v,$$ 
satisfies the relation 
$$Div\, X_{u_1,\ldots,u_n}=$$
$$\sum_{i=1}^n(-1)^{n+i-1}\der_i
(id\wedge\der_1\wedge \cdots\wedge\hat{\der_i}\wedge \cdots \wedge \der_n)
(u_1,\ldots,u_{n}))=$$
$$
\sum_{i=1}^n(-1)^{n+i-1}(\der_i\wedge\der_1\wedge \cdots\wedge \hat{\der_i}\wedge \cdots \wedge \der_n)
(u_1,\ldots,u_{n}))=$$
$$
(-1)n\sum_{i=1}^n(\der_1\wedge \cdots\wedge\der_i\wedge \cdots \wedge \der_n)
(u_1,\ldots,u_{n}))=$$
$$(-1)^nn\,Jac_n^S(u_1,\ldots,u_n).$$
In other words,
$$Div\,X_{u_1,\ldots,u_{n-1}}+n\,R_{u_1,\ldots,u_{n-1}}=0.$$ 
So, any interior derivation is an element of the subspace 
$\{X-{n}^{-1}Div\,X: X\in W_n(U)\}.$ 

\begin{lm} \label{10aug}  
Let $U$ be an associative commutative algebra with 
the space of commuting derivations ${\D}=\{\der_1,\ldots,\der_n\}.$ 
Then for any $X\in W_n$ the linear map $g=g_X:U\rightarrow U$ given by 
$g(u)=X(u)+\lambda \,Div\,X\,u$ is a derivation of the 
algebra $(U,Jac_{n+1}^W),$ if and only if $\lambda=-n^{-1}.$  
\end{lm}

{\bf Proof.} Consider the polynomial 
$$f(X,u_1,\ldots,u_{n+1})=$$ 
$$g_X(Jac_{n+1}^W(u_1,\ldots,u_{n+1}))-\sum_{i=1}^{n+1}
Jac_{n+1}^W(u_1,\ldots,u_{i-1},g_X(u_i),u_{i+1},\ldots,u_n).$$
We should prove that  $f(u_1,\ldots,u_{n+1})=0$ for any $u_1,\ldots,
u_{n+1}\in U,$ if and only if $\lambda=-n^{-1}  .$ 

Notice that $g=g_X$ is $\D$-invariant for any $X\in W_n$: 
$$(\rho(\der_i)g)(u)=\der_i(g(u))-g(\der_i(u))=$$
$$\der_i(Xu)+\lambda \der_i(Div\, X\, u)-X(\der_i(u))-
\lambda\, Div\,X\,\der_i(u)=$$
$$\der_i(X)u+\lambda\, \der_i(X)\, u=g(\der_iX)u.$$

We know that $Jac_{n+1}^S$ is $\D$-invariant: 
$$\rho(\der_i)(id\wedge \der_1\wedge\cdots\wedge\der_n)=0.$$
Hence by corollary \ref{22sept}
$f$ as a linear combination of compositions of $\D$-invariant 
polynomials, is also $\D$-invariant. 

As we mentioned before $Jac_{n+1}^W$ is graded. 
Notice  that $g$ has grade degree $|X|$: 
$$g(U_{[s]})\subseteq U_{[s]+|X|}.$$
Therefore, the polynomial $f$ is also graded. 

So, we can use theorem~\ref{altynbala} to prove that $f=0.$

Notice that $f$ is skew-symmetric in $n$ parameters $u_1,\ldots,u_{n+1}.$ 
Therefore, for homogeneous $X\in W_n$ and 
$u_1,\ldots,u_{n+1}\in U,$ 
$$f(X,u_1,\ldots,u_{n+1})\in U_{[0]},$$ 
if and only if
$$|X|=-1, |u_1|+\ldots+|u_{n+1}|=0$$
or
$$|X|=0, |u_l|=-1 \mbox{\;\;for some $l$ and\;\;}
|u_1|+\cdots\hat{|u_l|}+\cdots+ |u_{n+1}|=0$$

In other words, if 
$$(x^{\alpha(0)}\der_i,x^{\alpha(1)},\ldots,x^{\alpha(n+1)})\in Support(f),$$
then 
$$|\alpha(0)|=0, |\alpha(1)|+\cdots+ |\alpha(n+1)|=n+1$$
or 
$$|\alpha(0)|=1, |\alpha(l)|=0,
\mbox{\; for some $l$ and} \;
 |\alpha(1)|+\cdots+\widehat{|\alpha(l)|}+\cdots
+|\alpha(n+1)|=n$$

So, we should check that 
$$f(\der_i,x^{\alpha(1)},\ldots,x^{\alpha(n+1)})=0,$$
if $|\alpha(1)|+\cdots+|\alpha(n+1)|=n,$
or
$$f(x_j\der_i,x^{\alpha(1)},\ldots,x^{\alpha(n+1)})=0,$$
if $\alpha(l)=0,$ for some $l$ and 
$$|\alpha(1)|+\cdots
+\hat{|\alpha(l)|}+\cdots+
+|\alpha(n~+1)|=n, |\alpha(s)|>0, s\ne l.$$
The first statement is evident: since $\der_i$ is a derivation of $U,$ 
$Div\,\der_i=0,$  
and $[\der_i,\der_j]=0$ for all $i$ and $j,$ then 
$f(\der_i,x^{\alpha(1)},\ldots,x^{\alpha(n+1)})=0,$ 
for any $\alpha(1),\ldots,\alpha(n+1).$ 
Since $f$ is skew symmetric in $u_i$ type parameters, in the second case  
we can take, say,  $l=1.$ Then $|\alpha(2)|=\ldots=|\alpha(n+1)|=1.$ 
So, in the second case we have  
$$f(x_j\der_i,1,x_1,\ldots,x_n)=$$

$$x_j\der_i(Jac_{n+1}^W(1,x_1,\ldots,x_n))
+\lambda\,\delta_{i,j}(Jac_{n+1}^W(1,x_1,\ldots,x_n))$$
$$-Jac_{n+1}^W(\lambda\delta_{i,j} 1,x_1,\ldots,x_n)$$
$$
-\sum_{s=1}^{n}
Jac_{n+1}^W(1,x_1,\ldots,x_{s-1},x_j\der_i(x_s)+\lambda\, \delta_{i,j}x_s,
x_{s+1},\ldots,x_n)=$$

$$x_j\der_i(Jac_{n}^S(x_1,\ldots,x_n))
+\lambda\,\delta_{i,j}(Jac_{n}^S(x_1,\ldots,x_n))$$
$$-\lambda\,\delta_{i,j}Jac_{n}^S(x_1,\ldots,x_n)$$
$$
-\sum_{s=1}^{n}
Jac_{n}^S(x_1,\ldots,x_{s-1},x_j\der_i(x_s)+\lambda\, \delta_{i,j}x_s,
x_{s+1},\ldots,x_n)=$$

$$
-\sum_{s=1}^{n}\delta_{i,s}
Jac_{n}^S(x_1,\ldots,x_{s-1},x_j,x_{s+1},\ldots,x_n)$$
$$
-\sum_{s=1}^n\lambda\, \delta_{i,j}Jac_n^S(x_1,\ldots,x_{s-1},x_s,
x_{s+1},\ldots,x_n)=$$

$$-Jac_n^S(x_1,\ldots,x_{i-1},x_j,x_{i+1},\ldots,x_n)
-\lambda\,n\delta_{i,j}=$$

$$-\delta_{i,j}(1+\lambda n).$$
Thus $f=0,$ iff $\lambda=-n^{-1}.$ 

\begin{crl} \label{14aug} 
Any divergenceless derivation of $U$ is also a derivation 
of $(U,Jac_{n+1}^W)$.    
\end{crl}

\begin{lm} \label{der_iD} 
Suppose that $D(1)=0, D(x_i)=0$ for any $i=1,\ldots,n.$ Then 
$D(\der_i(u))=\der_i(D(u))$ for any $u\in U.$
\end{lm}

{\bf Proof.} Since 
$$Jac_{n+1}^W(u,1,x_1,\ldots,\hat{x_i},\ldots,x_n)=\der_i(u),$$
then, by the Leibniz rule
$$D(\der_i(u))=
D(Jac_{n+1}^W(u,1,x_1,\ldots,\hat{x_i},\ldots,x_n))=$$
$$Jac_{n+1}^W(D(u),1,x_1,\ldots,\hat{x_i},\ldots,x_n)=\der_i(D(u)),$$
for any $u\in U.$  

\begin{lm} \label{10augu} 
Let $U=K[x_1^{\pm 1},\ldots,x_n^{\pm 1}].$ 
Suppose that $D\in Der\,(U,Jac_{n+1}^W)$ and $D(u)=0,$ for any $u\in U^+.$ 
Then $D(u)=0,$ for any $u\in U.$
\end{lm}

{\bf Proof.} 
Suppose that $D(x^{\alpha})\ne 0$ for some $\alpha\in 
{\bf Z}^n\setminus {\bf Z}^n_+$ and $D(x^{\beta})=0,$ for all $\beta\in {\bf Z}^n, |\beta|>|\alpha|.$ If $\alpha_i+\epsilon_i<0,$ for some $i,$ 
then we obtain a contradiction: by lemma \ref{der_iD} in this case 
$$|\alpha+\epsilon_i|>|\alpha|\Rightarrow D(x^{\alpha+\epsilon_i})=0
\Rightarrow 
D(x^{\alpha})=(\alpha_i+1)^{-1}D(\der_i(x^{\alpha+\epsilon_i}))=0.$$
Thus, $\alpha_i=-1,$ for all $i=1,\ldots,n.$ 

Notice that 
$$Jac_{n+1}^W(u,x_1^2,\ldots,x_n^2)=
2^{n-1} x^{\theta}(2-\sum_{i=1}^nx_i\der_i)u,$$
for any $u\in U.$
In particular,
$$Jac_{n+1}^W(x^{-2\theta},x_1^2,\ldots,x_n^2)= 2^n (n+1)x^{-\theta}.$$
Thus,
$$D(x^{-\theta})=$$
$$2^{-n}(n+1)^{-1}D(Jac_{n+1}^W(x^{-2\theta},x_1^2,\ldots,x_n^2)=$$
$$2^{-n}(n+1)^{-1}Jac_{n+1}^W(D(x^{-2\theta}),x_1^2,\ldots,x_n^2)=$$
$$2^{-1}(n+1)^{-1}x^{\theta}(2-\sum_{i=1}^nx_i\der_i)D(x^{-2\theta}).$$
On the other hand,
$$|x_i^{-1}|> |x^{-\theta}|\Rightarrow D(x^{-1}_i)=0, i=1,\ldots,n,$$
$$Jac_{n+1}^W(1,x_1^{-1},\ldots,x_n^{-1})=(-1)^nx^{-2\theta},$$
and by the Leibniz rule 
$$D(x^{-2\theta})=$$
$$(-1)^nD(Jac_{n+1}^W(1,x_1^{-1},\ldots,x_n^{-1}))=0.$$
Therefore, $D(x^{-\theta})=0.$ Contradiction. 

\begin{lm} \label{10augus} 
Any derivation of the algebra $(U,Jac_{n+1}^W)$ has the form 
$X-n^{-1}Div\,X$ for some $X\in W_n(U).$ 
\end{lm}

{\bf Proof.} By lemma \ref{10augu} any derivation of $F\in 
Der\,(U,Jac_{n+1}^W)$ is defined by its restriction on $U^+=K[x_1^{\pm 1},\ldots,x_n^{\pm 1}].$ 
Let $D$ be the restriction of derivation $F$ on $U^+.$ 

{\bf Step 1}. 
Prove that 
\begin{equation} \label{13a}
D(1)\in <x^{\alpha}: \alpha\in {\bf Z}^n: \alpha\ne -\theta>.
\end{equation}

Notice that $|x^{-\theta}|=-2n\ne 0.$ Therefore, if 
$\sum_{i=1}^nx_i\der_iD(1)= 0,$ then (\ref{13a}) is true. 

Consider now the case $\sum_{i=1}^nx_i\der_iD(1)\ne 0.$
Since $1=Jac_{n+1}^W(1,x_1,\ldots,x_n),$ by the Leibniz rule we have 
$$D(1)=$$
$$Jac_{n+1}^W(D(1),x_1,\ldots,x_n)+\sum_{i=1}^n
Jac_{n+1}^W(1,x_1,\ldots,x_{i-1},D(x_i),x_{i+1},\ldots,x_n)=$$
$$(1-\sum_{i=1}^nx_i\der_i)D(1)+\sum_{i=1}^n
Jac_{n}^S(x_1,\ldots,x_{i-1},D(x_i),x_{i+1},\ldots,x_n).$$
Thus, 
$$\sum_{i=1}^n
Jac_{n}^S(x_1,\ldots,x_{i-1},D(x_i),x_{i+1},\ldots,x_n)=
\sum_{i=1}^nx_i\der_i D(1).$$
Therefore, 
$D(1)$ is a linear combination of elements of the form 
$Jac_n^S(u_1,\ldots,u_n)$ and  
by lemma \ref{3augu},  (\ref{13a}) is true. 

So, according to (\ref{13a})
$R=\sum_{\alpha\in {\bf Z}^n, \alpha\ne -\theta}
\lambda_{\alpha}x^{\alpha},$
for some $\lambda_{\alpha}\in K.$ 
If $\lambda_{\alpha}\ne 0,$ then there exists some $i=i(\alpha),$ such 
that $\alpha_i\ne -1.$ Therefore,   by lemma \ref{10aug} there exists  
a derivation of $(U,Jac_{n+1}^W)$ of the form 
 $G:=\sum_{\alpha\in {\bf Z}^n}
\lambda_{\alpha}(n(\alpha_i+1)^{-1}x^{\alpha+\epsilon_i}\der_i-
x^{\alpha}).$ Take $D'=D+G.$
It is easy to see that $$D'\in Der\,(U,Jac_{n+1}^W)$$ and $$D'(1)=0.$$ 
  
{\bf Step 2.}
Suppose that $D(x^{\bar\alpha})\ne 0,$  for 
some $\bar\alpha\in {\bf Z}^n_+, |\bar\alpha|=t,$ but 
$D(x^{\beta}=0,$ for all $\beta\in{\bf Z}^n_+,$ such that $|\beta|<t.$  
As we have shown above (step 1) we can assume that $t>0.$
Consider separately the cases $t=1$ and $t>1.$ 

{\bf The case $t=1.$} 
By the Leibniz rule for $D\in Der\,(U,Jac_{n+1}^W),$ we have 
$$D(Jac_n^S(u_1,\ldots,u_n))-
\sum_{l=1}^nJac_n^S(u_1,\ldots,u_{l-1},D(u_l),u_{l+1},\ldots,u_n)=$$
$$D(Jac_{n+1}^W(1,u_1,\ldots,u_n))-
\sum_{l=1}^nJac_{n+1}^W(1,u_1,\ldots,u_{l-1},D(u_l),u_{l+1},\ldots,u_n)=$$
$$Jac_{n+1}^W(D(1),u_1,\ldots,u_n).$$
Thus $D\in Der\,(U,Jac_n^S),$ if $D(1)=0.$ By theorem \ref{S}, 
$D$ as a derivation of $U$ under multiplication $Jac_n^S$ 
can be presented in the form 
$$D=\bar D+D_0+ \tilde D,$$
where 
$$\bar D=
\lambda_{-\theta}D_{-\theta},$$
$$D_0=\lambda_0\Delta,$$
$$\tilde D=\sum_{i=1}^n\lambda_ix^{-\theta+\epsilon_i}\der_i+\sum_{i<j, \alpha\in {\bf Z}^n}
\lambda_{i,j}(\alpha)D_{i,j}(x^{\alpha}),$$
for some $\lambda_0,\lambda_{-\theta}, \lambda_i, \lambda_{i,j}.$ 

Notice that $\tilde D$ is a divergence-free differential operator of 
first order. By corollary \ref{14aug} any divergence-free derivation  $D$ is a 
derivation of $(U, Jac_{n+1}^W).$ Since $D(1)=0,$ $\lambda_0=0.$  
So, $\bar D=D-\tilde D$ is also a derivation of $(U,Jac_{n+1}^W).$ 
 
Prove that $\bar D=0.$ As we mentioned before,
$$Jac_{n+1}^W(x^{-\theta},x_1,\ldots,x_n)=(n+1) x^{-\theta}.$$
Therefore,
\begin{equation}\label{13august}
\bar D(x^{-\theta})=(n+1)^{-1}\bar D(Jac_{n+1}^W(x^{-\theta},x_1,\ldots,x_n)).
\end{equation}
Let $\pi: U\rightarrow <1>$ be projection to $<1>.$ 
Notice that 
$$\pi(\bar D(x^{-\theta}))=\lambda_{-\theta},$$
$$\pi(Jac_{n+1}^W(\bar D(x^{-\theta}),x_1,\ldots,x_n))=\lambda_{-\theta},$$
Take projections to $<1>$ from the both parts of (\ref{13august}).  
We find that 
$$\lambda_{-\theta}=(n+1)^{-1}\lambda_{-\theta}.$$
So, $n\lambda_{-\theta}=0,$ and $\lambda_{-\theta}=0,$ $\bar D=0.$  

{\bf The case $t>1.$}  Prove that this case is not possible. 
By lemma~\ref{der_iD} $\der_iD(x^{\bar\alpha})=0,$ 
for any $i=1,\ldots,n.$ So, $D(x^{\bar\alpha})=\lambda \cdot 1,$ for some 
$0\ne \lambda\in K.$ 

Notice that 
$$Jac_{n+1}^W(u,x_1,\ldots,x_n)=(1-\sum_{i=1}^nx_i\der_i)u,$$
for any $u\in U.$ 
Therefore, by the Leibniz rule, 
$$(1-|\bar\alpha|)\lambda=(1-|\bar\alpha|)D(x^{\bar\alpha})=$$
$$D(Jac_{n+1}^W(x^{\bar\alpha},x_1,\ldots,x_n))=$$
$$Jac_{n+1}^W(D(x^{\bar\alpha}),x_1,\ldots,x_n)=\lambda.$$
So, $t\lambda=0,$ and $\lambda=0.$ Contradiction.

Summarize the results obtained in all these steps and cases. 
We see that any derivation of $(U,Jac_{n+1}^W)$ can be 
presented as a linear combination of derivations of the form $X-n^{-1}Div\,X,$ 
where $X\in W_n(U).$ 

{\bf Proof of theorem \ref{W}}. Suppose that $J$ is an ideal of $U$ under 
$(n+1)$-multiplication $\omega=Jac_{n+1}^W.$ Then $J$ is ideal of $U$ under 
$n$-multiplication $i(1)\omega=Jac_n^S.$ Since $(U,Jac_n^S)$ is simple \cite{Filipov1}, $J=0$ or $J=U.$ So, $(U,Jac_{n+1}^W)$ is simple. 

Other statements of theorem \ref{W} follow from lemma \ref{10augus}.


\begin{thebibliography}{10}

\bibitem{Cautheron} P. Cautheron, {\em Some remarcs concerning Nambu 
mechanics,} Lett. Math. Phys., {\bf 37}(1996), 103-116. 

\bibitem{DalTakh} Y.L. Daletskii, L.A. Takhtajan, {\em Leibniz and Lie algebra
structures for Nambu algebra,} Lett. Math. Phys., {\bf 39}(1997),
127-141.

\bibitem{Dzhumavestnik} A.S. Dzhumadil'daev,
{\em A remark on spaces of invariant differential operators}, Vestnik Moskov.Univ, Ser1, mat.,mech.,  1982, No.2, 49-54 = engl.transl. Moscow Univ. Math. Bull., {\bf 37}(1982), No.2, p.63-68.

\bibitem{DzhCentral} A.S. Dzhumadil'daev,
{\em Central extensions of infinite-dimensional Lie algebras},  
Funct.Anal.Appl., {\bf 26}(1992), No.4, p.21-29 = engl.transl. 246-253.

\bibitem{Filipov1} V.T. Filippov, {\em $n$-Lie algebras,} Sib. Mat.
J. {\bf 26}(1985), No.6, 126-140 = engl.transl. 
Siberian Math.J., {\bf 26}(1985), No.6, 879-891.

\bibitem{Filipov2} V.T. Filippov, {\em On $n$-Lie algebra of
jacobians,} Sib.Mat.J., {\bf 39}(1998), No.3, 660-669 = engl.transl. 
Siberian Math. J., {\bf 39}(1998), No.3, 573-581.

\bibitem{Kurosh} A.G. Kurosh, {\em Multi-operator rings and
algebras,} Uspechi Matem.Nauk,  {\bf 24}(1969), No.1, 3-15. 

\bibitem{Nambu} Y.Nambu, {\em Generalized Hamiltonian mechanics,}
Phys. Rev., D 7, 2405-2412, 1973.

\bibitem{Takh} L.A. Takhtajan, {\em On the foundation of the
generalized Nambu mechanics,} Commun. Math. Phys., {\bf
160}(1994), 295-315.

\bibitem{Takh1} L.A. Takhtajan, {\em A higher order ananog of the Chevalley-Eilenberg complex and the deformation theory of $n$-gebres,} Algebra i Analis, {\bf 61}(1994), No.2, 262-272 = engl.transl. St. Petersburg Math. J., {\bf 6}(1995), No.2, 429-438.
\end{thebibliography}
\end{document}